\magnification=1200
\input amstex
\catcode`\@=11
\ifx\amspptloaded@AmS\relax\catcode`\@=\active
 \endinput\else\let\amspptloaded@AmS\relax\fi
\parindent10\p@
\ifnum\mag=1000
  \hsize 16.5 cm
  \vsize 23 cm
  \captionwidth@\hsize
  \advance\captionwidth@-3 cm
\else  
  \hsize 16.5  true cm
  \vsize 23 true cm
  \captionwidth@\hsize
  \advance\captionwidth@-3 true cm
\fi
\normallineskiplimit\p@
\font@\ninerm=cmr9
\font@\eightrm=cmr8
\font@\sixrm=cmr6
\font@\ninei=cmmi9    \skewchar\ninei='177
\font@\eighti=cmmi8   \skewchar\eighti='177
\font@\sixi=cmmi6     \skewchar\sixi='177
\font@\ninesy=cmsy9   \skewchar\ninesy='60
\font@\eightsy=cmsy8  \skewchar\eightsy='60
\font@\sixsy=cmsy6    \skewchar\sixsy='60
\font@\ninebf=cmbx9
\font@\eightbf=cmbx8
\font@\sixbf=cmbx6
\font@\nineit=cmti9
\font@\eightit=cmti8
\font@\ninesl=cmsl9
\font@\eightsl=cmsl8
\font@\tensmc=cmcsc10
\font@\lc=cmcsc10 
\def\tenpoint{\def\pointsize@{10}%
 \normalbaselineskip 13.333333 pt
 \abovedisplayskip 13.333333 pt plus 3.333333 pt minus 9.999999 pt
 \belowdisplayskip 13.333333 pt plus 3.333333 pt minus 9.999999 pt
 \abovedisplayshortskip 0 pt plus 3.333333 pt
 \belowdisplayshortskip 7.777777 pt plus 3.333333 pt minus 4.444444 pt
 \smallskipamount= 3.333333 pt plus 1.111111 pt minus 1.111111 pt
 \medskipamount=6.666666 pt plus 2.222222 pt minus 2.222222 pt
 \bigskipamount=13.333333 pt plus 4.444444 pt minus 4.444444 pt
 \parskip=2 pt plus 1 pt
 \parindent 16.666666 pt
 \textonlyfont@\rm\tenrm
 \textonlyfont@\it\tenit
 \textonlyfont@\sl\tensl
 \textonlyfont@\bf\tenbf
 \textonlyfont@\smc\tensmc
 \ifsyntax@\def\big##1{{\hbox{$\left##1\right.$}}}\else
 \let\big\tenbig@
 \textfont\z@\tenrm  \scriptfont\z@\sevenrm  \scriptscriptfont\z@\fiverm
 \textfont\@ne\teni  \scriptfont\@ne\seveni  \scriptscriptfont\@ne\fivei
 \textfont\tw@\tensy \scriptfont\tw@\sevensy \scriptscriptfont\tw@\fivesy
 \textfont\thr@@\tenex \scriptfont\thr@@\tenex \scriptscriptfont\thr@@\tenex
 \textfont\itfam\tenit
 \textfont\slfam\tensl
 \textfont\bffam\tenbf \scriptfont\bffam\sevenbf
  \scriptscriptfont\bffam\fivebf
 \fi
 \setbox\strutbox\hbox{\vrule height8.5\p@ depth3.5\p@ width\z@}%
 \setbox\strutbox@\hbox{\vrule height8\p@ depth3\p@ width\z@}%
 \normalbaselines\tenrm\ex@=.2326ex}
\def\eightpoint{\def\pointsize@{8}%
 \normalbaselineskip10\p@
 \abovedisplayskip10\p@ plus2.4\p@ minus7.2\p@
 \belowdisplayskip10\p@ plus2.4\p@ minus7.2\p@
 \abovedisplayshortskip\z@ plus2.4\p@
 \belowdisplayshortskip5.6\p@ plus2.4\p@ minus3.2\p@
 \smallskipamount=3\p@ plus 1\p@ minus 1\p@
 \medskipamount=6\p@ plus 2\p@ minus 2\p@
 \bigskipamount=12\p@ plus 4\p@ minus 4\p@
 \parskip=0\p@ plus 1\p@
 \parindent 11.111111 pt
  \textonlyfont@\rm\eightrm
 \textonlyfont@\it\eightit
 \textonlyfont@\sl\eightsl
 \textonlyfont@\bf\eightbf
 \ifsyntax@\def\big##1{{\hbox{$\left##1\right.$}}}\else
 \let\big\eightbig@
 \textfont\z@\eightrm \scriptfont\z@\sixrm  \scriptscriptfont\z@\fiverm
 \textfont\@ne\eighti \scriptfont\@ne\sixi  \scriptscriptfont\@ne\fivei
 \textfont\tw@\eightsy \scriptfont\tw@\sixsy \scriptscriptfont\tw@\fivesy
 \textfont\thr@@\tenex \scriptfont\thr@@\tenex \scriptscriptfont\thr@@\tenex
 \textfont\itfam\eightit
 \textfont\slfam\eightsl
 \textfont\bffam\eightbf \scriptfont\bffam\sixbf
   \scriptscriptfont\bffam\fivebf
 \fi
 \setbox\strutbox\hbox{\vrule height7\p@ depth3\p@ width\z@}%
 \setbox\strutbox@\hbox{\vrule height6.5\p@ depth2.5\p@ width\z@}%
 \normalbaselines\eightrm\ex@=.2326ex}
\def\tenbig@#1{{\hbox{$\left#1\vbox to8.5\p@{}\right.\n@space$}}}
\def\eightbig@#1{{\hbox{$\textfont\z@\ninerm\textfont\tw@\ninesy
 \left#1\vbox to6.5\p@{}\right.\n@space$}}}
\def\footmarkform@#1{$^{#1}$}
\let\thefootnotemark\footmarkform@
\def\makefootnote@#1#2{\insert\footins
 {\interlinepenalty\interfootnotelinepenalty
 \eightpoint\splittopskip\ht\strutbox\splitmaxdepth\dp\strutbox
 \floatingpenalty\@MM\leftskip\z@\rightskip\z@\spaceskip\z@\xspaceskip\z@
 \indent{#1}\footstrut\ignorespaces#2\unskip\lower\dp\strutbox
 \vbox to\dp\strutbox{}}}
\newcount\footmarkcount@
\footmarkcount@\z@
\def\footnotemark{\let\@sf\empty\relaxnext@\ifhmode\edef
 \@sf{\spacefactor\the\spacefactor}\/\fi
 \def\next@{\ifx[\next\let\next\nextii@\else
  \ifx"\next\let\next\nextiii@\else
  \let\next\nextiv@\fi\fi\next}%
 \def\nextii@[##1]{\footmarkform@{##1}\@sf}%
 \def\nextiii@"##1"{{##1}\@sf}%
 \def\nextiv@{\global\advance\footmarkcount@\@ne
  \footmarkform@{\number\footmarkcount@}\@sf}%
 \futurelet\next\next@}
\def\footnotetext{\relaxnext@
 \def\next@{\ifx[\next\let\next\nextii@\else
  \ifx"\next\let\next\nextiii@\else
  \let\next\nextiv@\fi\fi\next}%
 \def\nextii@[##1]##2{\makefootnote@{\footmarkform@{##1}}{##2}}%
 \def\nextiii@"##1"##2{\makefootnote@{##1}{##2}}%
 \def\nextiv@##1{\makefootnote@{\footmarkform@{\number\footmarkcount@}}{##1}}%
 \futurelet\next\next@}
\def\footnote{\let\@sf\empty\relaxnext@\ifhmode\edef
 \@sf{\spacefactor\the\spacefactor}\/\fi
 \def\next@{\ifx[\next\let\next\nextii@\else
  \ifx"\next\let\next\nextiii@\else
  \let\next\nextiv@\fi\fi\next}%
 \def\nextii@[##1]##2{\footnotemark[##1]\footnotetext[##1]{##2}}%
 \def\nextiii@"##1"##2{\footnotemark"##1"\footnotetext"##1"{##2}}%
 \def\nextiv@##1{\footnotemark\footnotetext{##1}}%
 \futurelet\next\next@}
\def\adjustfootnotemark#1{\advance\footmarkcount@#1\relax}
\let\topmatter\relax
\newbox\titlebox@
\setbox\titlebox@\vbox{}
\Invalid@\overlong
\def\overlong@{\def\next@{\ifx\next\overlong\def\filhss@
 {plus\@m\p@ minus\@m\p@}\def\next@\overlong{\nextii@}\else
 \def\filhss@{plus\@m\p@\relax}\let\next@\nextii@\fi\next@}}
\def\title{\relaxnext@
 \def\nextii@##1\endtitle{{\let\\=\cr
 \global\setbox\titlebox@\vbox{\tabskip\z@\filhss@
 \halign to\hsize{\tenpoint\bf\hfil\ignorespaces####\unskip\hfil\cr##1\cr}}}}%
 \overlong@
 \futurelet\next\next@}
\newif\ifauthor@
\newbox\authorbox@
\def\author{\relaxnext@
 \def\nextii@##1\endauthor{{\let\\=\cr
 \global\setbox\authorbox@\vbox{\tabskip\z@\filhss@
 \halign to\hsize{\tenpoint\smc\hfil\ignorespaces####\unskip\hfil\cr##1\cr
 }}}}\overlong@\global\author@true
 \futurelet\next\next@}
\newif\ifaffil@
\newbox\affilbox@
\def\affil{\relaxnext@
 \def\nextii@{\bgroup\let\\=\cr
 \global\setbox\affilbox@\vbox\bgroup\tabskip\z@\filhss@
 \halign to\hsize\bgroup\tenpoint\hfil\ignorespaces####\unskip\hfil\cr}%
 \overlong@
 \global\affil@true
 \futurelet\next\next@}
\def\endaffil{\cr\egroup\egroup\egroup}
\newcount\addresscount@
\addresscount@\@ne
\def\address#1{\expandafter\gdef\csname address\number\addresscount@
 \endcsname{\noindent\eightpoint\ignorespaces#1\par}\global
 \advance\addresscount@\@ne}
\newif\ifdate@
\def\date#1{\global\date@true\gdef\date@{\tenpoint\ignorespaces#1\unskip}}
\newif\ifthanks@
\def\thanks#1{\global\thanks@true\gdef\thanks@{\eightpoint\ignorespaces
 #1\unskip}}
\Invalid@\nofrills
\Invalid@\usualspace
\newif\ifnofrills@
\def\usualspace@#1{\ifnofrills@\def\usualspace{#1}\fi}
\def\nofrills@#1#2{\def\next@{\ifx\next\nofrills\nofrills@true\let#2\relax
 \def\next@\nofrills{\nextii@}\else\nofrills@false
 \def#2{#1}\let\next@\nextii@\fi\next@}}
\newif\ifkeywords@
\def\thekeywords@{}
\def\keywords{\relaxnext@\nofrills@{{\eightpoint\bf Key words:\enspace}\rm }%
\keywords@
 \def\nextii@##1{\def\thekeywords@{\usualspace@{{\it\enspace}}\noindent
  \eightpoint\keywords@\ignorespaces##1\par}}%
\global\keywords@true
\futurelet\next\next@}
\newif\ifsubjclass@
\def\thesubjclass@{}
\def\subjclass{\relaxnext@\nofrills@{{{\bf 1991 Mathematics subject
 classifications:\enspace}\rm }}\subjclass@
 \def\nextii@##1{\def\thesubjclass@{\usualspace@
  {{\rm\spacefactor2000 \space}}\noindent\eightpoint
  \subjclass@\ignorespaces##1\par}}%
 \global\subjclass@true
 \futurelet\next\next@}
\def\proclaim{\innerproclaim@}
\def\endproclaim{\innerendproclaim@}
\newif\ifabstract@
\def\theabstract@{}
\def\abstract{\relaxnext@\nofrills@{{\bf Abstract.\enspace}}\abstract@
 \long\def\nextii@##1{\long\gdef\theabstract@{\usualspace@
  {{\eightpoint\enspace}}\eightpoint\abstract@\ignorespaces##1\par}}%
 \global\abstract@true
 \futurelet\next\next@}
\def\pretitle{}
\def\preauthor{}
\def\preaffil{}
\def\predate{}
\def\preabstract{}
\def\prepaper{}
\def\endtopmatter{\hrule height\z@\vskip-\topskip
 \pretitle
 \vskip 26.666666 pt plus 13.333333 pt minus 13.333333 pt
 \unvbox\titlebox@
 \preauthor
 \ifauthor@\vskip 13.333333 pt plus 6.666666pt minus 3.333333 pt%
\unvbox\authorbox@\fi
 \preaffil
 \ifaffil@\vskip10\p@ plus5\p@ minus2\p@\unvbox\affilbox@\fi
 \predate
 \ifdate@\vskip6\p@ plus2\p@ minus\p@\hbox to\hsize{\hfil\date@\hfil}\fi
 \preabstract
 \ifthanks@\makefootnote@{}{\thanks@}\fi
 \ifabstract@\vskip 16.666666 pt plus 7
 pt minus 7
 pt
   {\eightpoint\indent  \theabstract@}\fi
 \ifsubjclass@\vskip 10
 pt plus 4.333333 pt minus 4.333333 pt
   {\eightpoint\indent\thesubjclass@}\fi
 \ifkeywords@\vskip 10
 pt plus 4.333333 pt minus 4.333333 pt
   {\eightpoint\indent\thekeywords@}\fi
   \prepaper
 \outer\def\proclaim{\innerproclaim@}%
 \outer\def\endproclaim{\innerendproclaim@}%
 \vskip 20 pt plus 13.333333 pt minus 6.666666 pt \tenpoint}
\newcount\addressnum@
\outer\def\enddocument{\nobreak\sfcode`\.=3000 \vskip12\p@ minus6\p@%
\nobreak\vskip12\p@ minus6\p@\addressnum@\z@
 \loop\ifnum\addressnum@<\addresscount@\advance\addressnum@\@ne
 \csname address\number\addressnum@\endcsname\repeat
 \vfill\supereject\end}
\newbox\headingbox@
\outer\def\heading{\relaxnext@
 \def\nextii@{\bigbreak\bgroup\let\\=\cr
 \global\setbox\headingbox@\vbox\bgroup\tabskip\z@\filhss@
 \halign to\hsize\bgroup\tenpoint\smc\hfil\ignorespaces####\unskip\hfil\cr}%
 \overlong@
 \futurelet\next\next@}
\def\endheading{\cr\egroup\egroup\egroup\unvbox\headingbox@
 \nobreak\medskip}
\def\subheading{\relaxnext@\nofrills@{.\enspace}\subheading@
 \def\nextii@##1{\medbreak\indent{\usualspace@{{\bf\enspace}}%
  \tenpoint\bf\ignorespaces##1\unskip\subheading@}\ignorespaces}%
 \futurelet\next\next@}
\newif\ifproclaim@
\def\innerproclaim@{\relaxnext@\nofrills@{.\enspace}\proclaim@
 \def\nextii@##1{\message{(##1)}\medbreak\indent\def\next{8}%
  \ifx\pointsize@\next\uppercase
  {\usualspace@{{\rm\enspace}}\ignorespaces\rm##1\unskip\proclaim@}\else
  {\usualspace@{{\smc\enspace}}\smc\ignorespaces##1\unskip\proclaim@}\fi\sl
  \ifproclaim@\Err@{Previous \expandafter
 \eat@\string\\proclaim has no matching \expandafter
 \eat@\string\\endproclaim}\else\proclaim@true\fi\ignorespaces}%
 \futurelet\next\next@}
\def\innerendproclaim@{\proclaim@false\par\rm
 \ifdim\lastskip<\medskipamount\removelastskip\penalty55 \medskip\fi}
\def\demo{\relaxnext@\nofrills@{.\enspace}\demo@
 \def\nextii@##1{\par\ifdim\lastskip<\smallskipamount\removelastskip
  \smallskip\fi\indent{\usualspace@{{\sl\enspace}}%
  \sl\ignorespaces##1\unskip\demo@}\rm
  \ifproclaim@\Err@{Previous \expandafter
  \eat@\string\\proclaim had no matching \expandafter
  \eat@\string\\endproclaim}\fi\ignorespaces}%
 \futurelet\next\next@}

\def\qed{\ifhmode\unskip\nobreak\fi\ifmmode\ifinner\else\hskip5\p@\fi\fi
 \hbox{\hskip5\p@\lower 1.5 pt\hbox{\vrule width .2 pt 
\vbox{\hrule width 4 pt height .2 pt \vskip 7.1 pt\hrule width 4 pt 
height .2 pt }\unskip\vrule width .2 pt}\hskip\p@}}
\def\cite#1{\relaxnext@
 \def\nextiii@##1,##2\end@{[{\bf##1},##2]}%
 \in@,{#1}\ifin@\def\next{\nextiii@#1\end@}\else
 \def\next{[{\bf#1}]}\fi\next}
\newcount\rostercount@
\newif\iffirstitem@
\newtoks\everypartoks@
\let\plainitem@\item
\def\par@{\everypartoks@=\expandafter{\the\everypar}\everypar{}}
\def\roster{\edef\leftskip@{\leftskip\the\leftskip}\relaxnext@
 \rostercount@\z@\def\item{\futurelet\next\rosteritem@}%
 \def\next@{\ifx\next\runinitem\let\next\nextii@\else
  \let\next\nextiii@\fi\next}%
 \def\nextii@\runinitem{\unskip
  \def\next@{\ifx\next[\let\next\nextii@\else
   \ifx\next"\let\next\nextiii@\else\let\next\nextiv@\fi\fi\next}%
  \def\nextii@[####1]{\rostercount@####1\relax
   \enspace{\rm(\number\rostercount@)}~\ignorespaces}%
  \def\nextiii@"####1"{\enspace{\rm####1}~\ignorespaces}%
  \def\nextiv@{\enspace{\rm(1)}\rostercount@\@ne~}%
  \par@\firstitem@false
  \futurelet\next\next@}%
 \def\nextiii@{\par\par@\penalty\@m\smallskip 
  \firstitem@true}%
 \futurelet\next\next@}
\def\rosteritem@{\iffirstitem@\firstitem@false\else\par\fi
 \leftskip3\parindent\noindent
 \def\next@[##1]{\rostercount@##1\relax
  \llap{\hbox to2.5\parindent{\hss\rm(\romannumeral\rostercount@)}\hskip
  .5\parindent}\ignorespaces}%
 \def\nextii@"##1"{%
  \llap{\hbox to2.5\parindent{\hss\rm##1}\hskip.5\parindent}\ignorespaces}%
 \def\nextiii@{\advance\rostercount@\@ne
  \llap{\hbox to2.5\parindent{\hss\rm(\romannumeral\rostercount@)}\hskip
  .5\parindent}}%
 \ifx\next[\let\next\next@\else\ifx\next"\let\next\nextii@\else
 \let\next\nextiii@\fi\fi\next}

\newif\ifnextRunin@
\def\endroster{\relaxnext@\par\leftskip@
 \penalty-50 \vskip-\parskip\smallskip
 \def\next@{\ifx\next\Runinitem\let\next@\relax
  \else\nextRunin@false\let\item\plainitem@\ifx\next\par
  \def\next@\par{\everypar=\expandafter{\the\everypartoks@}}%
  \else\def\next@{\noindent\everypar=\expandafter{\the\everypartoks@}}%
  \fi\fi\next@}%
 \futurelet\next\next@}
\newcount\rosterhangafter@
\def\Runinitem#1\roster\runinitem{\relaxnext@\rostercount@\z@
 \def\item{\futurelet\next\rosteritem@}%
 \def\runinitem@{#1}%
 \def\next@{\ifx\next[\let\next\nextii@\else\ifx\next"\let\next\nextiii@
  \else\let\next\nextiv@\fi\fi\next}%
 \def\nextii@[##1]{\rostercount@##1\relax\def\item@{{\rm(\romannumeral
  \rostercount@)}}\nextv@}%
 \def\nextiii@"##1"{\def\item@{{\rm##1}}\nextv@}%
 \def\nextiv@{\advance\rostercount@\@ne\def\item@{{\rm(\romannumeral
  \rostercount@)}}\nextv@}%
 \def\nextv@{\setbox\z@\vbox
  {\ifnextRunin@\noindent\fi
  \runinitem@\unskip\enspace\item@~\par
  \global\rosterhangafter@\prevgraf}%
  \firstitem@false\ifnextRunin@\else\par\fi
  \hangafter\rosterhangafter@\hangindent3\parindent
  \ifnextRunin@\noindent\fi\runinitem@\unskip\enspace
  \item@~\ifnextRunin@\else\par@\fi\nextRunin@true\ignorespaces}%
 \futurelet\next\next@}
\outer\def\Refs{\relaxnext@\def\refskip@{\hskip\@ne sp\hskip\m@ne sp}%
 \def\next@{\ifx\next\nofrills\def\next@\nofrills{\nextii@}\else
  \def\next@{\nextii@{REFERENCES}}\fi\next@}%
 \def\nextii@##1{\bigbreak\hbox to\hsize{\hfil\eightpoint
  \rm\ignorespaces##1\unskip\hfil}\nobreak
  \bigskip\eightpoint\sfcode`.=\@m}%
 \futurelet\next\next@}
\outer\def\Lit{\relaxnext@\def\refskip@{\hskip\@ne sp\hskip\m@ne sp}%
 \def\next@{\ifx\next\nofrills\def\next@\nofrills{\nextii@}\else
  \def\next@{\nextii@{LITERATURA}}\fi\next@}%
 \def\nextii@##1{\bigbreak\hbox to\hsize{\hfil\eightpoint
  \rm\ignorespaces##1\unskip\hfil}\nobreak
  \bigskip\eightpoint\sfcode`.=\@m}%
 \futurelet\next\next@}
\newbox\nobox@        \newbox\keybox@        \newbox\bybox@
\newbox\bysamebox@    \newbox\paperbox@      \newbox\paperinfobox@
\newbox\jourbox@      \newbox\volbox@        \newbox\issuebox@
\newbox\yrbox@                               \newbox\pagesbox@
\newbox\bookbox@
\newbox\bookinfobox@  \newbox\publbox@       \newbox\publaddrbox@
\newbox\finalinfobox@
\newif\ifno@          \newif\ifkey@          \newif\ifby@ \newif\ifmanyby@
\newif\ifbysame@      \newif\ifpaper@        \newif\ifpaperinfo@
\newif\ifjour@        \newif\ifvol@          \newif\ifissue@
\newif\ifyr@ \newif\iftoappear@              \newif\ifpages@ \newif\ifpage@
\newif\ifbook@ \newif\ifinbook@
\newif\ifbookinfo@    \newif\ifpubl@         \newif\ifpubladdr@
\newif\iffinalinfo@   \newif\ifafterbook@
\newif\iffirstref@    \newif\iflastref@      \newif\ifprevjour@
\newif\ifprevbook@    \newif\ifprevinbook@   \newif\ifnojourinfo@
\newdimen\maxbysamerule@
\maxbysamerule@1in
\def\ref@{\global\no@false\global\key@false\global\by@false
 \global\bysame@false\global\paper@false\global\paperinfo@false
 \global\jour@false\global\vol@false\global\issue@false
 \global\yr@false\global\toappear@false\global\pages@false\global\page@false
 \global\book@false\global\inbook@false
 \global\bookinfo@false\global\publ@false\global\publaddr@false
 \global\finalinfo@false
 \bgroup\ignorespaces}
\Invalid@\moreref
\outer\def\ref{\begingroup
 \noindent\hangindent20\p@\hangafter\@ne\firstref@true
 \lastref@false\def\moreref{\egroup\endref@\global\firstref@false\ref@}\ref@}
\def\refdef@#1#2{\def#1{\egroup
 \csname\expandafter\eat@\string#1@true\endcsname
 \expandafter\setbox
 \csname\expandafter\eat@\string#1box@\endcsname\hbox\bgroup#2}}
\refdef@\no\relax \refdef@\key\relax
\def\manyby{\egroup\global\manyby@true\by@true\setbox\bybox@\hbox\bgroup\smc}
\def\by{\egroup\by@true\bysame@false\global\manyby@false
 \setbox\bybox@\hbox\bgroup\rm}
\def\bysame{\egroup\bysame@true\bgroup}
\refdef@\paper\sl
\refdef@\paperinfo\relax
\def\jour{\egroup\jour@true\prevjour@true\setbox
 \jourbox@\hbox\bgroup}
\refdef@\vol\bf
\refdef@\issue\relax \refdef@\yr\relax
\def\toappear{\egroup\toappear@true\bgroup}
\refdef@\pages\relax
\def\page{\egroup\page@true\setbox\pagesbox@\hbox\bgroup}
\refdef@\book\relax
\def\inbook{\egroup\inbook@true\previnbook@true\setbox
 \bookbox@\hbox\bgroup}
\refdef@\bookinfo\relax
\refdef@\publ\relax
\refdef@\publaddr\relax
\refdef@\finalinfo\relax
\def\setpunct@{\def\prepunct@{\ifnum\lastpenalty<0
 \edef\penalty@{\penalty\the\lastpenalty}\unpenalty,\ifafterbook@''\fi
  \penalty@\relax\space\else
 \ifdim\lastskip=\@ne sp\unskip\unskip
 \edef\penalty@{\penalty\the\lastpenalty}\unpenalty,\ifafterbook@''\fi
  \penalty@\relax\space
 \else,\ifafterbook@''\fi\space\fi\fi\afterbook@false}}
\def\ppunbox@#1{\prepunct@\unhbox#1\unskip}
\def\endref@{\let\prepunct@\relax
 \iffirstref@
  \ifno@\hbox to 22.222222 pt {\hss[\unhbox\nobox@\unskip]\kern 8.5 pt }\else
  \ifkey@\hbox to 22.222222 pt {\hss\unhbox\keybox@\unskip\kern 8.5 pt} 
  \else \hbox to 8.5 pt{}\fi\fi
  \ifmanyby@
   \ifby@\hbox{\unhcopy\bybox@\unskip}\setpunct@
  \global\setbox\bysamebox@\hbox{\unhcopy\bybox@\unskip}\else
  \ifbysame@\ifdim\wd\bysamebox@>\maxbysamerule@
    \hbox to\maxbysamerule@{\leaders\hrule\hfill}\else
    \hbox to \wd\bysamebox@{\leaders\hrule\hfill}\fi\setpunct@\fi
   \fi
  \else
  \ifby@\unhcopy\bybox@\unskip\setpunct@\fi\fi
 \fi
 \ifpaper@\ppunbox@\paperbox@\setpunct@\fi
 \ifpaperinfo@\ppunbox@\paperinfobox@\setpunct@\fi
 \ifjour@\ppunbox@\jourbox@\setpunct@
   \ifvol@\ppunbox@\volbox@\unskip\setpunct@\fi
   \ifissue@\ \unhbox\issuebox@\unskip\setpunct@\fi
   \ifyr@\ (\unhbox\yrbox@\unskip)\setpunct@\fi
   \iftoappear@\ (to appear)\setpunct@\fi
   \ifpages@\prepunct@ pp.\ \unhbox\pagesbox@\unskip\setpunct@\fi
   \ifpage@\prepunct@ p.\ \unhbox\pagesbox@\unskip\setpunct@\fi
 \else
  \ifprevjour@\unskip\nojourinfo@false
   \ifvol@\else\ifissue@\else\ifyr@\else\nojourinfo@true\fi\fi\fi
   \ifnojourinfo@\else,\fi
   \ifvol@\ppunbox@\volbox@\unskip\setpunct@\fi
   \ifissue@\ \unhbox\issuebox@\unskip\setpunct@\fi
   \ifyr@\ (\unhbox\yrbox@\unskip)\setpunct@\fi
   \iftoappear@\ (to appear)\setpunct@\fi
   \ifpages@\prepunct@  pp.\ \unhbox\pagesbox@\unskip\setpunct@\fi
   \ifpage@\prepunct@ p.\ \unhbox\pagesbox@\unskip\setpunct@\fi
  \fi
 \fi
 \ifbook@\prepunct@``\unhbox\bookbox@\unskip\afterbook@true\setpunct@\fi
 \ifinbook@\prepunct@\unskip\ in ``\unhbox\bookbox@\unskip\afterbook@true
  \setpunct@\global\book@true\fi
 \ifbookinfo@\ppunbox@\bookinfobox@\setpunct@\fi
 \ifpubl@\ppunbox@\publbox@\setpunct@\fi
 \ifpubladdr@\ppunbox@\publaddrbox@\setpunct@\fi
 \ifbook@
  \ifyr@\prepunct@\unhbox\yrbox@\unskip\setpunct@\fi
  \iftoappear@\ifafterbook@''\fi\ (to appear)\afterbook@false
   \setpunct@\fi
  \ifpages@\prepunct@ pp.\ \unhbox\pagesbox@\unskip\setpunct@\fi
  \ifpage@\prepunct@ p.\ \unhbox\pagesbox@\unskip\setpunct@\fi
 \else
  \ifprevinbook@\unskip
   \ifyr@\prepunct@\unhbox\yrbox@\unskip\setpunct@\fi
   \iftoappear@\ (to appear)\setpunct@\fi
   \ifpages@\prepunct@ pp.\ \unhbox\pagesbox@\unskip\setpunct@\fi
   \ifpage@\prepunct@ p.\ \unhbox\pagesbox@\unskip\setpunct@\fi
  \fi
 \fi
\iffinalinfo@.\ifafterbook@''\fi\afterbook@false
\spacefactor3000\relax\space\unhbox\finalinfobox@\else
 \iflastref@.\ifafterbook@''\fi\afterbook@false\else;\ifafterbook@''\fi
  \afterbook@false\space\fi
\fi}
\def\endref{\egroup\global\lastref@true\endref@\global\prevjour@false\global
 \previnbook@false\par\endgroup}
\newif\iflogo@
\def\nologo{\logo@false}
\def\logo{\logo@true}
\logo@false
\newif\ifdraft@
\def\nodraft{\draft@false}
\def\draft{\draft@true}
\draft@false
\output={\output@}
\def\output@{%
\ifdraft@  \shipout\vbox{\vbox to\vsize
  {\boxmaxdepth\maxdepth\pagecontents}\baselineskip2pc
\hbox to\hsize{\relaxnext@\hphantom{\eightpoint\jobname\quad
 \number\month/\number\day/\number\year\quad DRAFT}  \hfil
  \tenpoint\number\pageno\hfil \eightpoint\jobname\quad
 \number\month/\number\day/\number\year\quad DRAFT
 }}\else\ifnum\pageno=\@ne\shipout\vbox{\vbox to\vsize
  {\boxmaxdepth\maxdepth\pagecontents}\baselineskip2pc
 \iflogo@\hbox to\hsize{\hfil\eightpoint Typeset by \AmSTeX}\fi}\else
 \shipout\vbox{\vbox to\vsize
  {\boxmaxdepth\maxdepth\pagecontents}\baselineskip2pc
 \hbox to\hsize{\hfil\tenpoint\number\pageno\hfil}}%
 \fi\fi
  \global\advance\pageno\@ne
 \ifnum\outputpenalty>-\@MM\else\dosupereject\fi}
 \def\footertext{{\hbox to\hsize{\hfil\tenpoint\number\pageno\hfil}}}
\def\footertext{\hbox to\hsize{\relaxnext@\hphantom{\eightpoint
\jobname\quad
  \number\month/\number\day/\number\year\quad DRAFT}\hfil
  \tenpoint\number\pageno\hfil
  \eightpoint\jobname\quad\number\month/\number\day/\number\year\quad DRAFT}}
\def\footnoterule{\vskip-3\p@\hrule width 2truein \vskip 2.6\p@}
\def\m@g{\mag\count@
  \hsize 6.5 true in\vsize 9.0 true in\dimen\footins 8 true in%
\captionwidth@\hsize
  \advance\captionwidth@-1.5 true in
}
\tenpoint
\catcode`\@=\active

\pageno=1
%
%

\define\super{^}
%
\define\E{\operatorname{E}}
\define\Var{\operatorname{Var}}
\redefine\P{\operatorname{Prob}}
\define\Th{\operatorname{Th}}
\redefine\o{\omega}

\redefine\b{\beta}

\define\ep{\epsilon}
\redefine\i{\infty}
\redefine\o{\omega}
\redefine\al{\alpha}

\redefine\S{\Sigma}
\redefine\i{\infty}

\redefine\df#1#2{\dsize\frac{#1}{#2}}
\redefine\db#1#2{\dsize\binom{#1}{#2}}

\define\){\right)}
\define\({\left(}
\define\plus{\oplus}
\redefine\a{\overline{a}}
\redefine\bb{\overline{b}}
\define\f{\overline{f}}
\redefine\p{\overline{p}}
\define\limin{{\liminf}_{n\to\infty}\P(n,\p;\psi)}
\define\limsu{{\limsup}_{n\to\infty}\P(n,\p;\psi)}

\define\LL{L_+}
\define\LLL{L_{\le}}
\define\cG{\Cal G}
\define\bcG{\bar{\Cal G}}
\define\pplus{\bar\plus}
\topmatter
\title
 CONVERGENCE IN HOMOGENEOUS  RANDOM GRAPHS
\endtitle
\author
 TOMASZ \L UCZAK\footnote"\dag"{\ Department
of Discrete Mathematics, Adam Mickiewicz University,
Pozna\'n, Poland. This work was written while the author visited
The Hebrew University of Jerusalem.}
and SAHARON SHELAH\footnote"\ddag"{\ Institute of Mathematics,
The Hebrew University of Jerusalem, Jerusalem, Israel.
Partially supported by the United States and Israel Binational
Science Foundation, publ. 435.}
\endauthor
\abstract{For a sequence $\p=(p(1),p(2),\dots) $ let $G(n,\p)$
denote the
random graph with vertex set $\{1,2,\dots,n\}$ in which two
vertices $i$, $j$
are adjacent with probability $p(|i-j|)$,  
independently for each pair.
We study  how the convergence of probabilities of  first order
properties of
$G(n,\p)$, can be affected by the behaviour of $\p$ and
  the strength of the language we use.}
\endtopmatter
\document
\nologo
\pageno=1
\NoBlackBoxes

\subheading{1. Introduction}

Random graph theory studies how probabilities of
properties of random graphs change when the size of the problem,
typically
the number of vertices of the random graph, approaches infinity.
The most commonly used random graph model is $G(n,p)$ the graph
with  vertex set $[n] =\{1,2,\dots,n\}$,
in which two vertices are joined by an
edge independently with probability $p$.
 It was shown  by Glebskii, Kogan, Liogonkii and Talanov
[GKLT 69] and, independently, by Fagin [Fa 76],
that in  $G(n,p)$,  the probability of every property
which can be expressed by a first order sentence $\psi$
tends to 0 or 1 as $n\to\infty$.
Lynch [Ly ~80]  proved that even if we add to the language
the successor  predicate   the probability of each first order
sentence still converges to a limit.  (Here and below
the probability of a sentence $\psi$ means the probability that
 $\psi$  is satisfied.)   However it is no longer true when we enrich
the language further. Kaufmann and Shelah [KS 85]
 showed the existence of a monadic second order
sentence $\phi$, which uses only the relation of adjacency in
$G(n,p)$, whose probability does not converge as $n\to\i$. Furthermore,
Compton, Henson and Shelah [CHS 87]
gave an example of a  first order sentence
$\psi$ containing predicate ``$\le$'' such that the probability of $\psi$
does not converge -- in fact in both these cases the probability of sentences
$\phi$
and $\psi$  approaches both 0 and 1  infinitely many times.

  One may ask whether analogous results remain valid when
 the probability $p_{ij}$ that  vertices
$i$ and $j$ are connected by an edge  varies with $i$ and $j$.
Is it still true that a zero-one law holds for every first order
property which uses only the adjacency relation?
Or, maybe, is it possible to put some restrictions on the set
$\{p_{ij} : i,j\in [n] \}$ such that the
convergence of the probability of each first order sentence is
preserved even in the case of a linear order?
The purpose of this paper is to shed some light  on  
problems of this type 
 in a model  of the random graph  in which the probability that two
vertices are adjacent may depend on their distance.

For a sequence  $\p=\{p(i)\}_{i=1}^{\i}$, where $0\le p(i)\le
1$, let $G(n,\p)$  be a graph with  vertex set
$[n] $, in which a pair of vertices $v,w\in [n] $
appears as an edge with probability $p(|v-w|)$,
independently for each pair (this is a finite  version of the
probabilistic model introduced by
Grimmett, Keane and Marstrand [GKM 84]).
  Furthermore, let us define $L$ as the
first order logic whose vocabulary contains the binary predicate of
adjacency relation, whereas in language $\LL$ the successor predicate is also
available and in $\LLL$ one may say that $x\le y$.
We study how the behaviour of  sequence $\p$ could affect the
convergence
of sequence $\{\P(n,\p;\psi)\}_{n=1}^{\i}$, where $\psi$ is a
sentence  from $L$, $\LL$ or $\LLL$  and
 $\P(n,\p;\psi)$ denotes  the probability that $\psi$ is satisfied
in a model with universe $[n]$,
adjacency relation determined by $G(n,\p)$, and,
in the case of languages $\LL$ and $\LLL$, additional
binary predicates ``$x=y+1$" and ``$x\le y$"
(here and below $x,y\in [n] $ are treated as natural numbers).
                  
The structure of the paper goes as follows. We start with
 the short list of basic notions and  results useful in the study
of first order theories.
  Then, in the next three sections, we study the convergence 
of sequence  $\P(n,\p;\psi)$, where $\psi$ is a 
first order sentence from languages $L$, $\LL$ and $\LLL$ 
respectively. It turns out that differences between those 
three languages are quite significant. 
Our first result gives a sufficient and necessary condition
which, imposed on $\p$, assures convergence of $\P(n,\p;\psi)$ 
 for every $\psi$ from $L$. In particular, we show 
that each sequence $\a$ can be ``diluted'' by adding  some 
additional zero terms in such a way that 
for the resulting sequence $\p$ and every $\psi$ 
from $L$ the probability  $\P(n,\p;\psi)$ tends either to 0 or to 1. 
It is no longer true for sentences $\psi$ from $\LL$. 
In this case  the condition $\prod_{i=1}^{\i} (1-p(i))>0$  
turns out to be sufficient (and, in a way, necessary) 
 for convergence of $\P(n,\p;\psi)$ for every $\psi$ from $\LL$.
Thus, the convergence of  $\P(n,\p;\psi)$ depends mainly on the 
positive terms of $\p$ and adding zeros to $\p$, in principle,  
does not improve convergence properties of $G(n,\p)$.
On the contrary, we give an example of a property $\psi$ from 
$\LL$ and a sequence $\a$ such that for every $\p$ obtained from $\a$ 
by adding enough zeros $\limin=0$ whereas $\limsu=1$. 

The fact that, unlike in the case of $L$,  additional zeros in $\p$
might spoil the convergence properties of $G(n,\p)$ becomes even more 
evident in the language  $\LLL$. We  show that there exists 
a property $\psi$ from $\LLL$ such that every infinite
sequence of positive numbers $\a$ can be diluted by adding 
zeros such that for the resulting sequence $\p$ we have
$\limin=0$ but $\limsu=1$. 
Furthermore, it turns out that in this case we can distinguish 
two types of the limit behaviour of $\P(n,\p;\psi )$.
If   $\prod_{i=1}^{\i} (1-p(i))^i>0$  then for  every $\psi$
from $\LLL$ the probability $\P(n,\p;\psi )$ converges. However, 
if we assume only that  $\prod_{i=1}^{\i} (1-p(i))>0$, 
then, although, as we mentioned above,  we can not assure the 
convergence of  $\P(n,\p;\psi )$ for every $\psi $ from $\LLL$, 
a kind of a ``weak convergence'' takes place, namely, for every
$\psi $ from $\LLL$ we have $\limsu-\limin<1$.

 We also  study convergence properties of a similar
random graph model $C(n,\p)$ which  uses as the universe a circuit
of $n$ points instead of interval $[n]$.
It appears that in this case convergence does not depend very much
on the strength of the language.
We show that there is a first order sentence $\psi$ which uses only
adjacency relation such that  for  every infinite sequence $\a$ of 
positive numbers, there exists a sequence $\p$, 
obtained from $\a$ by adding 
enough  zero terms, for which the probability that $C(n, \p)$ 
has $\psi$ does not converge.
On the other hand,  under some rather natural constraints
imposed on $\p$, one can  show that a zero-one law holds
for a large class of first order properties
provided   $\p$ is such that
$\log(\prod_{i=1}^n (1-p(i))/\log n\to 0$.

We  conclude  the paper with  some additional remarks
concerning presented results and their  possible  generalizations.
Here we show also that  even in the case when a zero-one law holds  
 the probability $\P(n,\p;\psi)$  can tend to the limit very slowly 
and no decision procedure can 
determine the limit $\lim_{n\to\i}\P(n,\p;\psi)$ 
for every first order sentence ~$\psi$.

\subheading{2. First order logic -- useful tools and basic facts}

In this part of the paper
we gather basic notions and facts concerning the first
order logic which shall be used later on.
Throughout this section ${\bar L}$  denotes
a first order  logic whose vocabulary contains a finite number
of  predicates $P_1,P_2,\dots,P_m$, where the
$i$-th predicate $P_i$ has $j_i$ arguments.
More formally, for a vocabulary $\tau$  let $L^\tau$  be
the first order logic (i.e the set of first order formulas)
in the vocabulary $\tau$. A  $\tau$-model  $M$, 
 called  also  a $L^\tau$-model,  is defined in the 
ordinary way. In the paper we  use
four vocabularies:
\roster
\item"1)" $\tau_0$  such that the $\tau_0$-models 
 are just graphs; we write $L$ instead $L^{\tau_0}$.
\item"2)"  $\tau_1$  such that  the  $\tau_1$-models
 are,  up to isomorphism,  
quintuples $([n], S,c,d, R) $, where $[n]=\{1,2\dots,n\}$, 
$S$ is the successor relation,
$c$, an individual constant, is 1, the other 
 individual constant $d$ is  $n$ and 
$([n],R)$ is a graph; we write $L_+$ instead $L^{\tau_1}$.
\item"3)"  $\tau_2$  such that the  $\tau_2$-models
 are, up to isomorphism, triples $([n],\le, R)$, where $\le $
is the usual order on $[n]$  and $([n],R)$ is a graph;  
we write $L_\le$ instead $L^{\tau_2}$.
\item"4)"  $\tau_3$  such that the  $\tau_3$-models
 are,  up to isomorphism, triples $([n],C, R)$,   
where $C$ is the ternary relation of between
in clockwise order (i.e. $C(v_1,v_2,v_3)$ means that 
$v_{\sigma(1)}\le v_{\sigma(2)}\le v_{\sigma(3)}$ for  
some cyclic permutation $\sigma $ of set $\{1,2,3\}$) 
and $([n],R)$ is a graph; 
 we write $L^c_\le$  instead $L^{\tau_3}$.
\endroster
For every natural number $k$ and  ${\bar L}$-model
$M=(U^M;P^M_1,P^M_2,\dots,P_m^M)$  we
set
$$\Th_k(M)=\{\phi:M\models\phi,\,\phi \text{ is a first order
sentence from 
${\bar L}$ of quantifier depth $\le k$}\}.$$

The {\sl Ehrenfeucht game} of length $k$ on two ${\bar L}$-models
${M^1}$ and ${M^2}$
is the game between two players, where
 in the $i$-th step of the game,  $i=1,2,\dots,k$, the first
player chooses a point $v^1_i$ from
$U^{{M^1}}$ or $v^2_i$ from $U^{{M\super 2}}$
and the second player  must answer by picking a point from the
universe of the other model.  The second player wins such a game
 when the structures induced by points
$v^1_1,v^1_2,\dots,v^1_k$ and
$v^2_1,v^2_2,\dots,v^2_k$ are isomorphic, i.e.
$$P^{{M^1}}_i(v^1_{l_1},v^1_{l_2},\dots,v^1_{l_{j_i}})=
 P^{{M^2}}_i(v^2_{l_1},v^2_{l_2},
\dots,v^2_{l_{j_i}})$$
for every $i=1,2,\dots,m$ and
$l_1,l_2,\dots, l_{j_i}\in [k] $.
The following well known fact
(see for example Gurevich [Gu 85])
 makes the Ehrenfeucht game a useful tool in studies
of first order properties.

\proclaim{Fact 1}
 Let  $M_1$, $M_2$  be $\tau$-models. 
Then the second player has a winning strategy for the Ehrenfeucht
game of length
$k$ played on ${M^1}$ and ${M^2}$ if and only if 
$\Th_k({M^1})=\Th_k({M^2})$.\qed
\endproclaim

In the paper we shall use also  a ``local'' version of the above fact.
For a ${\bar L}$-model and a point $v\in U^M$
the {\sl neighbourhood} $N(v)$ of $v$ is the set of all points
  $w$ from  $U^M$ such that either $v=w$ or 
there exists $i$,  $ i=1,2,\dots,m$,
and  $v_1, v_2,\dots, v_{j_i}\in U^M$  for which
$ P^M_i(v_1, v_2,\dots, v_{j_i})$,  where $v=v_{i_1}$,
 $w=v_{i_2}$ for some $i_1,i_2$.
Set $N_1(v)=N(v)$ and  $N_{i+1}(v)=\bigcup_{w\in N_i(v)}N(w)$ 
for $i=1,2,\dots$, and define  the distance between two points 
$v,w\in U^M$ as the smallest $i$ for which $v\in N_i(w)$.
Clearly, the distance defined in such a way is a symmetric function
for which the triangle inequality holds. 
Now let ${M^1}$ and ${M^2}$ be two ${\bar L}$-models,
$v^1\in U^{{M^1}}$ and $v^2\in U\super
{{M^2}}$.
We say that pair $({M^1},v^1)$ is {\sl $k$-equivalent}
to $({M^2},v^2)$ when the second
player has a winning strategy in the ``restricted'' Ehrenfeucht game
of length $k$ in which, in the first step players must choose 
vertices $v^1=v^1_1$ and $v^2=v^2_2$
and in the $i$-th step, $i=1,2,\dots,k$ they are forced to
pick vertices $v^1_i$ and $v^2_i$ from sets
$\bigcup_{j=1}^{i-1}N_{3^{k-i}}(v^1_j)$
and  $\bigcup_{j=1}^{i-1}N_{3^{k-i}}(v^2_j)$.
The following result (which, in fact, is a version of a special case
of Gaifman's result from [Ga 82])  is an easy consequence of Fact 1.

\proclaim{Fact 2}
Let ${M^1}$ and ${M^2}$ be two ${\bar L}$-models
such that for  $l=1,2,\dots,k$,  $i=1,2$,
  every choice of points
$v^i_1, v^i_2,\dots,v^i_l\in U^{{M^i}}$
and $v^{3-i}_1,v^{3-i}_2,\dots,v^{3-i}_{l-1}
\in U^{{M^{3-i}}}$, such that  no two of the 
$v_j^i$ and $v_j^{3-i}$ are 
 at a distance less than  $3^{k-l+1}$
from each other and $(M^i,v^i_j)$ is $(k-l+1)$-equivalent
to $(M^{3-i}, v_j^{3-i})$ for $j=1,2,\dots,l-1$, 
there exists $v_l^{3-i}\in U^{M^{3-i}}$
such that $(M^i,v^i_l)$ is $(k-l)$-equivalent to 
$(M^{3-i},v^{3-i}_l)$.

Then  $\Th_k({M^1})=\Th_k({M^2})$.\qed
\endproclaim

Finally, we need some results from the theory of additivity of
models.
Call $\S$  a {\sl scheme of a generalized sum with respect to
vocabularies $\hat\tau$,  $\tau$ and $\tau'$}
if for each predicate $P(\bar x)$
of  $\tau'$ 
and breaking $<\hskip-0.11cm\bar x_i\hskip-0.15cm>_{i\le k}$
of $\bar x$,  $\S$  gives  a first order, quantifier free formula 
$\phi_P (z_1,...,z_k)$ in vocabulary
$$\hat\tau\cup\{R_i^{\psi}\,:\,\text{$\psi$ is a
 quantifier free formula in   $L^\tau$ 
with free variables $\bar x_i$}\},$$
where $R_i^{\psi}$ denotes a zero-place predicate, i.e.
 a    truth value.  


\proclaim{Definition}
Let $\S$ be  a  scheme of a generalized sum with respect to
vocabularies $\hat\tau$,  $\tau$ and $\tau'$, 
$I$ be a $\hat\tau$-model and  $\{M_i\}_{i\in I}$
be a family of $\tau$-models. We shall say that a
  $\tau'$-model $N$
is a {\sl $(I,\S)$-sum of   $\{M_i\}_{i\in I}$} if 
the universe of $N$ is the disjoint sum of the universes of
$\{M_i\}_{i\in I}$ and for each
$\tau'$-predicate $P(\bar x) $ relation
$P^N$ is the set of $\bar a$
such that for some breaking $<\hskip-0.11cm\bar x_i\hskip-
0.15cm>_{i\le k}$
of $\bar x$  there are distinct members $t_1,t_2,\dots, t_k$ of $I$
and sequences $\bar a_i$ of members of $M_{t_i}$ of the same
length as $\bar x_i$,  $ i\le k$,  such that
$\bar a$ is a concatenation of $\bar a_1,\bar a_2,\dots,\bar a_k$
and, if for each $i\le k$ we interpret  $R^{\psi}_i$ as the truth
whenever $M_{t_i}\models \psi(\bar a_i)$,  in the model $I$  the formula
 $\phi_P   (t_1,t_2,...,t_k)$  is satisfied.
\endproclaim



\proclaim{Example 1}
\roster
\item"(i)"
Let $G_1, G_2, \dots, G_k$ be  graphs, treated as models of
language $L$,
whose vocabulary contains only a binary predicate interpreted as
the adjacency  relation. Then the graph
$$G=G_1\plus  G_2\plus \dots G_k\,,$$
defined as the sum of disjoint copies of $G_1,G_2,\dots, G_k$, is a
$(\S,I)$-sum of these graphs, for $I=\{1,2,\dots, k\}$
and  empty vocabulary $\hat\tau$.
\item"(ii)"
For  $i=1,2,\dots,m$, let $G_i$ be a graph with  vertex set
$\{1,2,\dots,n_i\}$  and
$$G=G_1\pplus  G_2\pplus \dots \pplus G_m$$
denote a graph with  vertex set $\{1,2, \dots,\sum_{i=1}\super
m n_i\}$
such that vertices ~$v$ and ~$w$ are adjacent in $G$ if and only if
for some $j=1,2,\dots,m$,
$$\sum_{i=1}^{j-1}n_i< v<w\le \sum_{i=1}^{j}n_i$$
and vertices $v-\sum_{i=1}^{j-1}n_i$ and $w-
\sum_{i=1}^{j-1}n_i$
are adjacent in $G_j$.

Let us view graphs $G_1,G_2, \dots,G_m,G$  as models  of
language $\LL$, which contains the adjacency relation and the
successor predicate and two individual constants which represent 
the first and the last elements of a graph.  
Then,  $G$ can be treated as a $(\S,I)$-sum of
$G_1,G_2,\dots, G_m$, for   $I=\{1,2,\dots,m\}$ and $\hat\tau=\LL$.
\item"(iii)"
Let graphs $G_1,G_2, \dots,G_m,G$  be defined as in the previous
case. Then, if these graphs are treated as models of 
language $\LLL$, which contains the 
adjacency relation and the predicate ``$\le$'', $G$ can be viewed
as a $(\S,I)$-sum of $G_1,G_2,\dots, G_m$, where
$I=\{1,2,\dots,m\}$ is the
model of  linear order.
\item"(iv)"
Let graphs $G_1,G_2, \dots,G_m,G$  be defined as in (ii) and
$\LLL^c$ be the language which contains predicate
$C(v_1,v_2,v_3)$ which
means that for some cyclic permutation $\sigma$ of indices 1,2,3
we have
$v_{\sigma(1)}\le v_{\sigma(2)}\le v_{\sigma(3)}$.
Then, if we treat $G_1, G_2,\dots, G_m$ as $\LLL$-models and $G$
as a $\LLL^c$-model,  $G$ can be viewed as $(\S,I)$-sum of
 $G_1, G_2,\dots, G_m$ with $I=\{1,2,\dots,m\}$
treated as a $\LLL^c$-model.
\endroster
\endproclaim

\demo{Remark}
Note that in  the definition of
a scheme of generalized sum the formula $\phi_P$  which
corresponds to predicate $P$ must be quantifier free.
This is the reason why we need  two individual constants
in  the language $\LL$.  

The main theorem  about $(\S,I)$-sums we shall use can be stated
as follows.

\proclaim{Fact 3} Let $\S$ be a fixed scheme of addition 
with respect to some  fixed vocabularies. 
Then, for every $k$ and
 $N$, a  $(\S,I)$-sum of $\{M_i\}_{i\in I}$,
 $\Th_k(N)$ can be computed
from $\{\Th_{k}(I,R_s)\,:\,s\in S\}$, where
$S=\{\Th_k(M)\,:\,M\text{ is a
$\tau$-model}\}$ and $R_s=\{i\in I\,:\, \Th_k(M_i)=s\}$.\qed
\endproclaim

We apply this result to $(\S,I)$-sums of graphs described  in
Example 1.

\proclaim{Fact 4} Let operations ``$\plus$''  and ``$\pplus$''
and languages $L$, $\LL$, $\LLL$ and $\LLL^c$ be defined in
the same way as in the Example 1. Furthermore, let $\cG$ and
$\bcG$ be families of graphs closed under $\plus$ 
and $\pplus$ respectively, and let $k$ be a natural number. 
  Then
\roster
\item"(i)"
there exists a graph $G\in \cG$ such that for every $H\in \cG$
$$\Th_k(G)=\Th_k(G\plus H)\,,$$
where in the above equation all graphs are treated as  $L$-models;
\item"(ii)"
 there exists a graph $\bar G\in \bcG$ such that
for every $\bar H\in \bcG$
$$\Th_k(\bar G)=\Th_k(\bar G\pplus \bar H\pplus\bar G)\,,$$
where either all graphs  are treated as $\LL$-models or all of them
are viewed as $\LLL$-models;
\item"(iii)"
there exists a graph $\bar G\in \bcG$ such that
for every $\bar H\in \bcG$
$$\Th_k(\bar G)=\Th_k(\bar G\pplus \bar H)\,,$$
where 
both $\bar G$ and $\bar G\pplus \bar H$ are treated as $\LLL^c$-models.
\endroster
\endproclaim

\demo{Proof} Let $\Cal U$ be  a set of all finite words over 
finite alphabet $S$. Words  from $\Cal U$
can be viewed as models of a language  $L(S)$, whose vocabulary
consists of   unary predicates ~$P_s$, for $s\in S$.  
Then, for any word  $\al$ which contains
$k$ copies of each letter of the alphabet and any other word $\b$
we have  $$\Th_k(\al)=\Th_k(\al\circ\b)\,,\tag 1$$
where $\al\circ\b$ denotes concatenation  of $\al$ and $\b$.

Let us set  $S=\{\Th_k(G): G\text{ is a $L$-model}\}$, and 
choose $\{G_s\}_{s\in S}$  such  that $\Th_k(G_s)=s$.
Furthermore, let 
$$G'= \bigoplus_{s\in S} G_s\quad\text{and}\quad 
G=\undersetbrace\text{$k$ times}\to{G'\plus\dots\plus G'}\;.$$
Then, from Fact 3 and (1), for every $H$, treated as 
a $L$-model, we have $\Th_k(G)=\Th_k(G\plus H)$. 

Now, treat  words from  $\Cal U$    as  
$\LL(S)$[$\LLL(S)$]-models for a language $\LL(S)$ [$\LLL(S)$] which, 
in addition to the predicates  $P_s$, contains also 
the successor predicate [the  predicate ``$\le$''].
It is not hard to show (see, for instance,  Shelah and Spencer [SS 94])
 that for every $k$ there exists a word $\al$ such that for every 
$\b$ we have $\Th_k(\al)=\Th_k(\al\circ\b\circ \al)$ for both 
$\LL(S)$ and $\LLL(S)$.
Thus, similarly as in the case of (i), 
the second part of the assertion follows from Fact~3.

Finally, let $\LLL^c(S)$ be a language which 
contains the  predicates $P_s$  and 
the ternary predicate   $C$ denoting clockwise order. 
 It is known (see again [SS 94])
that for every $k$ there exists a word $\al$ such that
for every other word $\b$ we have $\Th_k(\al)=\Th_k(\al\circ\b)$,
where this time both $\al$ and $\al\circ \b$ are treated 
as $\LLL^c(S)$-models. Hence, using  Fact~3 once again, 
we get the last part of Fact ~4.\qed

\subheading{3. Zero-one laws for language $ L$}

In this section we  characterize sequences $\p$ for which
the probability $\P(n,\p;\psi)$ converges for every sentence $\psi$
from $L$,
or, more precisely, for which a zero-one law holds,
i.e. it converges to either 0 or 1.
One could easily see that the proof of  either
Glebskii, Kogan, Liogonkii and 
Talanov [GKLT 69], or Fagin [Fa 76],  can be
mimicked whenever a  sequence $\p$  is such that
$$0<\liminf_{i\to\i}p(i)\le\limsup_{i\to\i}p(i)<1\,,$$
so it is enough to consider only the case when
$\liminf_{i\to\i}p(i)=0$ (if $\limsup _{i\to\i}p(i)=1$ 
one can instead consider properties of the complement of  $G(n,\p)$).
The main result of this section describes rather precisely how
 convergence properties of $G(n,\p)$ depend on the fact how fast
the product $\prod_{i=1}^{\i}(1-p(i))$ tends to 0.

\proclaim{Theorem 1}
\roster
\item"(i)" For every  sequence $\p=(p(1),p(2),\dots,)$ such that
$p(i)<1$ for all $i$ and
$$\log\Big(\prod_{i=1}^n (1-p(i))\Big)\big/\log n \to 0\tag 2$$
and every sentence $\psi $ from $L$  a zero-one law holds.
\item"(ii)" For every positive constant $\ep$ there exists a
sequence $\p$ and
 a sentence $\psi$ from $L$  such that
$$-\log\Big(\prod_{i=1}^n (1-p(i))\Big)\big/\log n <\ep$$
but $\limin=0 $  while $\limsu=1 $.
\endroster
\endproclaim

In order to show Theorem 1 we need some information about 
the structure of $G(n,\p)$.
A subgraph $H'$ of a graph $G$ with the vertex set
$[n] $ is the {\sl exact copy} of a graph $H$ with the vertex set
$[l]$, if for some $i$, $0\le i\le n-l$, and  every $j,k\in [l]$
the pair $\{j,k\}$ is an edge of $H$ if and only if $\{i+j,i+k\}$ 
appears as an edge of $H'$ and
$G$ contains no edges with precisely one end in 
$\{i+1,i+2,\dots,i+l\}$.
Furthermore, call a graph $H$ on the vertex set
$[l]$ {\sl admissible\/} by a sequence $\p$ if the probability that
$H=G(l,\p)$ is positive. We shall show first that,
 with probability tending to 1 as $n\to\i$, 
$G(n,\p)$ contains many disjoint exact
copies of every finite admissible graph, 
provided $\prod_{i=1}^n (1-p(i))$ tends to infinity slowly enough.

\proclaim{Lemma}
For $k\ge 1$ let  $\p$ be a sequence such that 
$$\prod_{i=1}^n (1-p(i))\ge n^{-1/(10k)},$$
and let $H$ be an admissible
graph with vertex set  $[k]$.  Then, 
  the probability that in  $G(n,\p)$ there exist  
at least $n^{0.1}$ vertex disjoint exact copies of $H$, 
none of them containing
vertices which are either less than $\log n$
or larger than $n-\log n$, tends to 1 as $n\to\i$.
\endproclaim

\demo{Proof}
Let $\Cal A$ denote the family of disjoint sets 
$A_i=\{\lceil\log n\rceil +ik+1,
\lceil\log n\rceil +ik+2, \dots,\lceil\log n\rceil +(i+1)k\}$, 
where $i=0,1,2, \dots,i_0-1,i_0=\lfloor (n-2\log n)/k\rfloor$.
For every set $A_i\in\Cal A$ the probability that the subgraph 
induced in
$G(n,\p)$ by $A_i$ is an exact copy of a graph $H$ with edge set 
$E(H)$ equals $P(A_i)=P(H)P'(A_i)$, where the factor
$$P(H)=\prod_{e=\{i,j\}\in E(H)}p(|i-j|)
\prod_{e'=\{i',j'\}\not\in E(H)}(1-p(|i'-j'|))> 0$$
remains the same for all sets $A_i$, whereas the probability
$P'(A_{i})$ that no vertices of $A_i$ are adjacent to 
vertices outside $A_i$, given by
$$P'(A_i)=\prod_{r=1}^k\prod_{s=1}^
{\lceil\log n\rceil+ik}(1-p(\lceil\log n\rceil+ik+r-s))
\prod_{t=\lceil\log n\rceil +(i+1)k+1}^n(1-p(t-\lceil\log n\rceil-ik-r))\,,$$
may vary with $i$. Nevertheless,  for a sequence $\p$ which fulfills
the assumptions of the Lemma, we have always
$$P'(A_i)\ge\Big( \prod_{r=1}^n (1-p(r))\Big)^{2k}\ge n^{-0.2}\;.$$
(Here and below we assume that all inequalities hold only for $n$ large
enough.) Thus, there exists a subfamily  $\Cal A'$ 
of $\Cal A$  with $\lfloor\sqrt n\rfloor$ elements,
such that for every $A\in \Cal A'$ the probability 
$P(A)$ is roughly the same, i.e. for   some function $f(n)$,
where $n^{-0.2}\le f(n)\le 1$,
$$f(n)(1-o(n^{-0.1}))\le P(A)\le f(n)(1+o(n^{-0.1}))$$
for all $A\in\Cal A'$.

Now, let $X$  denote the number of sets from $\Cal A'$ which are
exact copies of $H$. For  the expectation of $X$ we get
$$
\E X=\sum_{A\in \Cal A'}P(A)=(1+o(n^{-0.1}))f(n)
\lfloor\sqrt n\rfloor\ge n^{0.2}\,.$$
 To estimate the variance of $X$ we need to  find  an upper bound for
$$
\E X(X-1)=\sum_{\Sb A,B\in \Cal A'\\A\neq B\endSb}P(A)P(B)
\big/\(\prod_{r\in A}\prod_{s\in B}
 (1-p(|r-s|))\)\,.
$$
 Note first that $(1-n^{-0.2})^{n^{0.3}}\le n^{-1/(10k)}$, so in
every sequence $\p$ for which the assumption holds at most $n^{0.3}$ 
of the first $n$ terms are larger than $n^{-0.2}$. Hence,
for all, except for at most  $n^{0.8}$,  pairs $A,B\in \Cal A'$
we have 
$$\prod_{r\in A}\prod_{s\in B} (1-p(|r-s|))
\ge (1-n^{-0.2})^{k^2}\ge 1-n^{-0.1}\,.$$
Moreover,  for every $A,B\in \Cal A'$
$$\prod_{r\in A}\prod_{s\in B} (1-p(|r-s|))
\ge \prod_{i=1}^n (1-p(i))^2\ge n^{-1/(5k)}\ge n^{-0.1}\,.$$
Thus,
$$\E X(X-1)\le nf^2(n)(1+O(n^{-0.1}))
+n^{0.8}f^2(n)n^{0.1}\le nf^2(n)(1+O(n^{-0.1}))\,,$$
the variance of $X$ is $o((\E X)^2)$, and from
Chebyshev's inequality
with probability tending to 1 as $n\to\i$ we have $X>\E
X/2>n^{0.1}$.\qed

\demo{Proof of Theorem 1}
Let $\psi$ be a first order sentence  of quantifier depth $k$.
For two graphs $G_1$ and $G_2$ define graph $G_1\plus G_2$ as
the disjoint sum of $G_1$ and $G_2$.
Since all $p(i)$ are less than 1  the family of admissible graphs 
is closed under the operation ``$\plus$''. 
Thus, from Fact ~4, there exists an admissible  graph  $G$ such that
for every admissible graph $H$ we have
$\Th_k(G\plus H)=\Th_k(G)$
(Let us recall that all graphs are treated here as models of
language  $L$ which contains  one binary  predicate interpreted as
the adjacency relation.)
From the Lemma we know that, for every sequence  $\p$ for which  (2)
holds, the probability that $G(n,\p)$ contains an exact copy of $G$ 
tends to 1 as $n\to\i$. Thus, with probability $1-o(1)$,
$\Th_k(G(n,\p)) =\Th_k(G)$ and the first part of Theorem 1 follows.

Now fix  $k\ge \lceil 1/\ep \rceil$ and let $\bb$ be a sequence 
of natural numbers such that $b(1)> 6k$ and 
$b(m+1)\ge (b(m))^{50}$ (e.g. $b(m)=(2k)^{50^m}$). 
Let us define a sequence $\p$ setting 
$$p(i)=\cases 1/2& \quad\text{ for }\quad i\le b(1)\\
\frac{1}{3ik} &\quad\text{ for }\quad b(2m-1)< i\le b(2m), 
\quad \text{ where } m=1,2,\dots\\
0 &\quad\quad \text{otherwise}\,.
\endcases$$
Then,  using the fact that for every $x\in (0,2/3)$ 
$$\exp(-2x)< 1-x< \exp(-x)\;,$$
we get
$$\prod_{i=1}^n (1-p(i))\ge 2^{-b(1)} 
\prod_{i=1}^n\(1-\frac 1{3ki}\)\ge n^{-1/k}\ge n^{-\ep}\,.$$
 for every sufficiently large $n$. We shall show that 
the probability that $G(n,\p)$ contains an exact copy 
of the complete graph $K_l$
on $l=6k$ vertices approaches both 0 and 1 infinitely many times.
 
Indeed, for $n=b(2m+1)$ and $m$ large enough we have 
$$\align 
\prod_{i=1}^{n}(1-p(i))&\ge 2^{-b(1)} 
\prod_{i=1}^{b(2m)}\(1-\frac1{3ik}\)\ge O(1)(b(2m))^{-2/(3k)}\\
&\ge (b(2m+1))^{-1/(70k)}= n^{-1/(70k)}\;.\endalign 
$$ 
Thus, from the Lemma, the probability that $G(b(2m+1),\p)$ contains
an exact copy of $K_l$  tends to 1 as $m\to\i$. 
On the other hand,  the expected number of exact
copies of $K_l$ in $G(b(2m+2),\p)$ is, 
for $m$ large enough,  bounded from above by 
$$\align
b(2m+2)\Big(&\prod_{i=6k}^{b(2m+2)/2}(1-p(i))\Big)^{6k}
\le b(2m+2)\bigg(\prod_{i=b(2m+1)+1}^{b(2m+2)/2}
\(1-\df{1}{3ik}\)\bigg)^{6k}\\
&\le  (1+o(1)) b(2m+2)\(\df{b(2m+2)}{2b(2m+1)}\)^{-2}
\le (b(2m+2))^{-1/3}\;,
\endalign$$
and tends to 0 as $m\to\i$.\qed

 Note that  Theorem 1 implies that each sequence $\a$ can
be ``diluted'' by adding zeros so that for the
resulting sequence $\p$, graph $G(n,\p)$
 has  good convergence properties.

\proclaim{Corollary}
For every sequence $\a$, where $0\le a(i)<1$,
 there is a sequence $f(i)$ such that every sequence
$\p$ obtained from $\a$ by the addition of more than $f(i)$ zeros
after the $i$-th term of $\a$ and each sentence $\psi$ from $L$ 
a zero-one law holds. \qed
\endproclaim


\subheading{4. Convergence for language ${\LL}$}

In [Ly 80] Lynch showed that if  $p(i)$ does not depend on $i$,
i.e. when  $ p(i)=p_0$ for some constant 
$0<p_0<1$ and $i=1,2,\dots$, 
then $\P(n,\p;\psi)$ converges for every  sentence   $\psi$
from $\LL$. In fact,  his argument guarantees the existence of  
 $\lim_{n\to\infty}\P(n,\p;\psi)$ for each sentence $\psi$ from $\LL$
and every sequence $\p$ which tends to a positive limit strictly 
smaller than one. Furthermore, if 
$\liminf_{i\to\i}p(i)<\limsup_{i\to\i}p(i)$ the probability of the
property that vertices 1 and $n$ are adjacent does not converge,
so, it is enough to consider the case when $p(i)\to 0$. 
Our first result says that the condition  that
$\prod_{i=1}^{\i}(1-p(i))>0$, or, equivalently, 
$\sum_{i=1}^{\i}p(i)<\i$, is sufficient and,
in a way, necessary, for the convergence of $\P(n,\p;\psi)$ 
for every sentence $\psi$ from $\LL$.

\proclaim{Theorem 2}
\roster
\item"(i)" If $\p$ is a sequence such that $p(i)<1$ for all $i$  and
$$\prod_{i=1}^{\i}(1-p(i))>0\;,\tag3$$
then for every sentence $\psi$ from $\LL$
the limit $\lim_{n\to\i}\P(n,\p;\psi)$ exists.
\item"(ii)" 
For every function $\o(n)$ which tends to infinity as
$n\to\i$  there exist a sequence $\p$ and sentence $\psi$ 
from $\LL$ such that
$$\o(n)\prod_{i=1}^{n}(1-p(i))\to\i\tag4$$
but $\limin=0$ whereas $\limsu=1$.
\endroster
\endproclaim

\demo{Proof} We shall deduce the first part of Theorem 2
 from Fact 2. Since our language
contains the successor predicate, in this section the distance
between two vertices $v,w\in [n] $
of a graph $G$ will be defined as  the length of the shortest path
joining $v$  to $w$ in the graph $\hat G$ obtained from $G$ by 
adding to the set of
edges of $G$  all pairs $\{i,i+1\}$,  where $i=1,2,\dots,n-1$,
and the neighbourhood of a vertex $v$ will mean always
neighbourhood in ~$\hat G$.

Let $\p$ be a  sequence for which (3) holds and 
$\psi $ be a  sentence  from $\LL$ of quantifier depth $k$. 
Call a  pair $(H,v)$ {\sl safe}
if $H$ is an  admissible graph on  $[l]$ and 
$v$ is a vertex of $H$ which lies at a distance at least 
$3^k$ from both 1 and $l$. Since there are only finite
number of $k$-equivalence classes  we can find a finite 
family $\Cal H$ of safe pairs such that 
every safe pair $(H', v')$ is $k$-equivalent to some 
pair $(H,v)$ from $\Cal H$. Now, due to  the Lemma, 
with probability tending to 1 as $n\to\infty$, $G(n,\p)$ 
contains at least $k$ exact copies of every safe pair $(H,v)$ 
from $\Cal H$. Thus, roughly speaking, the ``local'' properties
of the ``internal'' vertices are roughly the same 
for all graphs $G(n,\p)$, provided $n$ is large enough. 

In order to deal with vertices lying near 1 and $n$ 
  we need to  ``classify'' 
graphs with vertex set $[n]$ according to their 
``boundary'' regions. 
  More specifically, 
let  $H_{n,k}(1)$ and $H_{n,k}(n)$ denote subgraphs induced in
$G(n,\p)$ by all vertices which lie within  the distance 
$3^{k+2}$ from 1 and $n$
respectively. We show that the probability that
$(H_{n,k}(1),1)$ [$(H_{n,k}(n),n)$]
is $(k+1)$-equivalent to  some $(H,v)$ converges as
$n\to\i$. 

Let $\ep$ be any positive constant.  Note first that the expected 
number of neighbours of a given vertex in $G(n,\p)$ 
is bounded from above by $C_1=2+\sum_{i=2}^\i p(i)$. 
Hence, the expected number of vertices in  $H_{n,k}(1)$ is 
less than $C_2=\sum_{i=0}^{3^{k+2}} C_1^i$ and, from Markov 
inequality, the probability that $H_{n,k}(1)$ contains more than 
$C_3=C_2/\ep$ vertices is less than $\ep$. 
Moreover, choose  $C_4$ in such a way that 
$\sum_{i\ge C_4}p(i)\le \ep/C_3$.
Then, the conditional probability that  some $v\in H_{n,k}(1)$ 
has a neighbour $w$ such that  $|v-w|\ge C_4$, provided the size 
of $H_{n,k}(1)$ is less than $C_3$,  is bounded from above by  $\ep$.
Hence, with probability at least $1-2\ep$, $H_{n,k}(1)$ 
 contains no vertices $v$ for which $v\ge C_5=3^{k+2}C_4+1$.
 Thus, for every $n,m\ge C_5$,
$H_{n,k}(1)$ and $H_{m,k}(1)$ are ``isomorphic'' with probability 
at least $1-4\ep$, or, more precisely, for every property $\phi$ 
$$|P(H_{n,k}(1)\text { has }\phi)-P(H_{m,k}(1)\text { has }\phi)|
\le 4\ep\;.$$
In particular, for every graph $H$ on $[l]$ vertices  
the probability that $(H_{n,k}(1),1)$ and $(H,1)$ are $(k+1)$-equivalent
converges as $n\to\i$. Clearly, the analogous result
 holds for $H_{n,k}(n)$.

To complete the proof note that the fact that $(H,v)$ 
and $(H',v')$ are $(k+1)$-equivalent implies that for every 
vertex $w$ in $(H,v)$,  lying within a distance 
$3^k$ from $v$,  there exists a vertex $w'$ in $(H',v')$ 
within a distance $3^k$ of $v'$ such that $(H,w)$ and $(H',w')$
are $k$-equivalent. Thus, if a graph $G$ with vertex set $[n]$ 
contains $k$ exact copies of every safe pairs from $\Cal H$,
$\Th_k(G)$ can be computed from the  $(k+1)$-equivalence 
classes of its $3^{k+1}$ neighbourhoods of $1$ and $n$
and  the assertion follows.

Now let $\o(n)$ be a function which tends to infinity  as $n\to\i$.
We may assume that $\o(n)$ is 
non-decreasing, and, say,  $\o(n)\le n/100$.
Let $f(2)=1$ and for $m\ge 3$ 
$$f(m)=\min\{l:\o(l)\ge 2^{f(m-1)}\}+ 4 m^3\;.$$
Define a sequence $\p$ setting
$$p(i)=\cases 1/m &\quad\text{ for }\quad 
f(m)-m^3\le  i\le f(m), \quad m=2,3,\dots \\
0&\quad\text{ otherwise}\,,\endcases $$
and let $m=\max\{l:f(l)\le n\}$. 
Note that for every $j$ we have $f(j)>f(j-1)+j^3$ so $p(i)$
is correctly defined. Furthermore, $f(j)\ge 2^{j-2}$ for all $j$, 
so  $\o(n)\ge 2^{2^{m-2}}$ and 
 the sequence $p(1),p(2), \dots, p(n)$ 
contains at most $m^4$ non-zero terms. 
Consequently, 
$$\o(n)\prod_{i=1}^n (1-p(i))\ge 2^{2^{m-2}}(1-1/2)^{m^4}
=2^{2^{m-2}-m^4}\to \i$$
as $m\to\i$. 
Furthermore, let $\psi$ be the property that vertices 1 and $n$ 
are joined by a path of length two. 
Then $P(2f(m)-2m^3,\p;\psi)=0$ while 
$$P(2f(m)-m^3,\p;\psi)\ge 1-\(1-1/{m^2}\)^{m^3-1}\ge 1-3e^{-m}
=1-o(1)\;.\qed$$

\demo{Remark} Note that in fact we have shown  that if
$\psi$ belongs to $\LL$ and a sequence $\p$ fulfills (3) then
$\lim_{n\to\i}\P(n,\p;\psi)$  is equal to the probability
that  $G(\i,\p)$ has $\psi$, where $G(\i,\p)$ is a graph 
with vertex set 
$V=V_1\cup V_2=\{1,2,\dots\}\cup\{\dots,-2,-1,\}$  which
contains no edges between sets $V_1$ and $V_2$ and vertices 
$v_i,w_i\in V_i$, $i=1,2$, are
adjacent with probability $p(|v_i-w_i|)$.

The second part of Theorem 2 suggests that the analog of the
Corollary derived from Theorem 1 is not valid for language $\LL$.
The following example shows that, in fact,
much more is true -- there exist
a sentence $\psi$ from $\LL$ and a sequence $\a$
such that for {\sl each} sequence $\p$ obtained from $\a$ by adding
enough zero terms the probability $\P(n,\p;\psi)$ 
approaches both 0 and 1 infinitely many times.

\proclaim{Example 2}
 Let $\bb$ be a sequence such that $b(j+1)\ge (b(j))^{10}$
(e.g. $b(j)=\exp(10^j)$) and 
$$a(i)=\cases i^{-0.2}&\quad\text{ for }\quad 
b(2j)< i \le b(2j+1)\\
 i^{-0.95}&\quad\text{ for }\quad b(2j+1)< i\le b(2j+2)\;.
\endcases$$ 
Furthermore,  let $\f$ be any sequence such that for every $i\ge 2$ we
have $f(i)>10\sum_{j=1}^{i-1}f(j)$,
$$p(j)=\cases a(i) &\quad\text{ if }\quad j=f(i)\\
0 &\quad\text{otherwise}\,.\endcases$$
and $\psi$ be a sentence from $\LL$ saying that
any two neighbours of vertex 1 are connected by a path of length
four not containing vertex 1.

Then $\limin=0$ and $\limsu=1$.
\endproclaim

\demo{Validation}
Let $n=f(m)$, where $m=b(2j+1)$, and let   $v=f(i)+1$ and $w=f(j)+1$ 
be  neighbours of 1 in $G(n,\p)$.
Note that $f(m)>4f(m-1)$, so for every edge $\{s,t\}$ of $G(n,\p)$ 
we have $|s-t|<n/4$, in particular   $v,w\le n/4$.
The probability that ~$v$ and ~$w$ are {\sl not} connected by a
path of type $v(v+f(k))(v+f(k)+f(j))(w+f(j))w$ is bounded from 
above by
$$\align
\prod_{k=1}^{m-1}&(1-a(k)a(j)a(k)a(j))\le
\prod_{r=m^{0.1}}^{m-1}(1-r^{-0.4}a^2(j))\\
&\le \exp\(-a^2(j)\sum_{r=m^{0.1}}^{m/2}
r^{-0.4}\)\le\exp\(-a^2(j)m^{0.5}\)\le \exp\(-m^{0.1}\)
\endalign$$
since for  a vertex  $u=f(l)$ either $l\ge m^{0.1}>b(2j)$ and 
then $a(l)= l^{-0.2}\ge m^{-0.2}$, or $l\le m^{0.1}$ and so
$a(l)\ge l^{-0.95}\ge m^{-0.2}$. 
In $G(n,\p)$ vertex 1 has at most  $m$ neighbours,
so  the expected number of pairs $v$ and $w$ such that both $v$ and
$w$ are adjacent to ~1 but they are not connected by a path of 
length 4 is bounded from above by $m^2\exp\(-m^{0.1}\)$ 
and tends to 0 as $m\to \infty$.

Now let  $n=f(m)$, where $m=b(2j+2)$, and let $v=f(i)+1$ and 
$w=f(j)+1$ denote the largest and the second largest neighbours
of 1 in $G(n,\p)$, respectively. Note that, since the sequence $\f$ 
grows very quickly,  every path $vv_1v_2v_3w$ of length four 
joining $v$ and $w$ has the property that among
$|v-v_1|$, $|v_1-v_2|$, $|v_2-v_3|$ and $|v_3-w|$ each from distances 
$f(i)$ and $f(j)$ appears once and some distance $f(k)$ 
appears twice. Since, for given 
$k$, there are at most eight possible paths of length four 
joining $v$ and $w$ whose ``edge lengths'' are $f(i), f(j), f(k), f(k)$,
the probability that $v$ and $w$ are joined by a path of length four 
is bounded from above by 
$$8\sum_{k=1}^m a(i)a(j)a(k) a(k)\;.$$ 
But, since $v$ and $w$ are largest neighbours of 1 in $G(n,\p)$,
with probability tending to 1 as $m\to\infty$ we have $i,j\ge m/2$
 and, consequently, $a(i),a(j)\le 2m^{-0.95}$. 
Thus, the probability that $v$ and $w$ are
connected by a path of length four is less than 
$$o(1)+32m^{-1.9}\sum_{k=1}^m (a(k))^2
\le o(1)+32m^{-1.9}\sum_{k=1}^m m^{-0.4}
\le o(1)+1/m=o(1)\;.\qed$$

In the last two sections we assumed  that $p(i)<1$ for all $i$. 
Now we discuss briefly  the situation when we allow the sequence 
$\p$ to contain terms which are equal to 1. 
If we are dealing with $\LL$, any finite number of ones in $\p$ 
is not a problem at all. Indeed, for each  formula $\psi$ one can 
easily find $\psi'$ such that $\psi $ holds in $G$ if and only if
$\psi'$ holds for $G'$,  where $G'$ is obtained from $G$ 
by deleting all edges $\{v,w\}$ for which 
$p(|v-w|)=1$.

When we use  language $L$ even a finite number of ones 
in $\p$ can cause some troubles. 
For example, if we set $p(1)=p(2)=1$ and the sequence $\p$ is such 
that for every $i< j<k$ we have  
$p(j-i)p(k-j)p(k-i)=0$ unless
$k=j+1=i+2$, then one can easily ``identify''  ~1 and $n$
in $G(n,\p)$  as those vertices which are contained in 
precisely one triangle. Thus, in this case, 
convergence properties for language  $L$ become 
very similar to those of $\LL$, in particular,
the assumption (2) in  Theorem 1i should be replaced by (3) (where 
products in (2) and (3) are taken over all $i$ such that  $p(i)<1$). 

When $\p$ contains an infinite number of ones 
the probability of some simple  properties of $G(n,\p)$ 
may  oscillate between 0  and 1.
To see it  set, for instance,
 $$p(i)=\cases 1 &\quad\text{ when } i=4^j \text{ and }j=0,1,\dots
   \\0 &\quad\text{ otherwise }\endcases$$
and let  $\phi$ be  the property that each edge of a graph is contained in
a cycle of length 4. Then, clearly,
$$\P(n,\p;\phi)=\cases 0  &\quad\text{ when } n=4^j+1
   \\1 &\quad\text{ otherwise }.\endcases$$                             

On the other hand we should mention that 
there are sequences with unbounded number of 
zeros and ones for which $\P(n,\p;\psi)$ converges.
Let us take for example the random sequence $\p_{\text{rand}}$ 
of zeros and ones such that 
$$P(p_{\text{rand}}(i)=0)= P(p_{\text{rand}}(i)=1)=1/2\;,$$ 
 independently for each $i=1,2,\dots$\,. 
Furthermore, for a given $k$, say that  
a graph $G$ has property  $\Cal A_k$ if for every subset $A$ 
of the vertices of $G$ with precisely $k$ elements and every 
$A'\subseteq A$ there exists a vertex $v$ of $G$ such that 
$v$ is adjacent to all vertices from $A'$ and not adjacent 
to all vertices from $A\setminus A'$. It is not hard to prove that
with probability 1 sequence $\p_{\text{rand}}$ 
has the property  that for every $k$ there exists $n_0=n_0(k)$ 
such that  $G(n,\p_{\text{rand}})$ has $\Cal A_k$ for $n\ge n_0$
(note that the probability space in this case is 
related only to the random construction  of  $\p_{\text{rand}}$ -- 
since this sequence  contains only zeros and ones 
once it is chosen graph $G(n,\p_{\text{rand}})$ is uniquely determined).
Thus,  the second player can easily win 
the Ehrenfeucht game of length $k$ on $G(n,\p_{\text{rand}} )$
and $G(m,\p_{\text{rand}} )$, provided  $n,m>n_0(k)$,
and so for every sentence $\psi$ from $L$ a zero-one law holds. 

\subheading{5.  Linear order case}

As one might expect, conditions which were sufficient 
for convergence of $\P(n,\p;\psi)$ for $\psi$ from $\LL$ 
are too weak to assure convergence for every $\psi$ from $\LLL$. 
Our first result states that, 
although for a sequence $\p$ with a finite number
of non-zero terms and every  $\psi$ from $\LLL$
the probability $\P(n,\p;\psi)$ tends to a limit as $n\to\i$,
the assumption of finiteness could not be replaced by any
convergence condition imposed on  positive terms 
of the sequence   $p(i)$.

\proclaim{Theorem 3}
\roster
\item"(i)" If $\p$ contains only finitely many  non-zero terms
then the probability $\P(n,\p;\psi)$ converges for every first order
sentence $\psi$ from $\LLL$.
\item"(ii)" For every infinite sequence $\a$, where
$a(i)>0$ for  $i=1,2,\dots$, and every positive constant $\ep>0$,
there exist a sequence $\p$ obtained from $\a$ by addition of some
zero terms  and a first order sentence $\psi$ from $\LLL$ such that
$\limsup_{n\to\i}\P(n,\p;\psi)=1$ and $ \liminf_{n\to\i}
\P(n,\p;\psi)\le \ep\,.$
\endroster
\endproclaim

\demo{Proof of Theorem 3}
Let $\p$ be a sequence with finitely many non-zero terms and let 
$\p'$ be obtained from $\p$ by replacing all terms equal to 
one by zeros. Since the successor relation can be expressed in $\LLL$,
for every sentence $\psi$ from $\LLL$ there exists 
 a sentence $\psi'$ in $\LLL$ such that $\P(n,\p;\psi)=\P(n,\p';\psi')$.
 Thus, we may assume that all terms of $\p$ are
strictly less than one. 

Now, let $\psi $ be a sentence of $\LLL$ of quantifier depth $k$
and let $C=C(k)>3^{k+1}$ be a constant such that no two 
vertices $v$, $v'$  of $G(n,\p)$ with $|v-v'|\le C-3^{k+1}$  
are, with positive probability,
joined by a path of length less than $3^{k+1}$ in $G(n,\p)$.
Now, in order to show the first part  of Theorem 3 it 
is enough to classify all graphs according to the structure of
the finite subgraphs induced by  subsets of vertices $\{v:v<C\}$
and $\{v: v>n-C\}$  
and observe that the Lemma implies that the probability 
that $G(n,\p)$ belongs to a given class converges as $n\to\i$.   
(Since Theorem 3i follows from much stronger Theorem~6 proven
below we omit details.)
  
Now let $\a$ be an infinite sequence of numbers such that 
$0<a(i)<1$ for $i=1,2,\dots$\,.  
Set  $f(1)=1$ and, for $i=2,3,\dots$,
$$ f(i)=\lceil\max\{(i+1)/a(i+1), 
4if(i-1) [1-\max\{a(j):j\le i-1\}]^{-(f(i-1))^2}\}\rceil\,.
\tag5$$
Moreover let
$$p(j)=\cases a(i)&\quad \text { if } j=f(i)\\
                0&\quad \text{ otherwise\,.}\endcases$$
Call a vertex $v$ of a graph a {\sl cutpoint} if a graph contains
no edges $\{w',w''\}$ such that $w'\le v$ and $w''> v$ and let
$\psi(r)$ be the property that a graph $G(n,\p)$ contains a 
cutpoint $v$ such that
 $$f(r)+r\le 2f(r)\le v\le n-2f(r)\le n-f(r)-r$$

Note first that the probability $\P(f(i),\p,\psi(r))$ tends
to 1 as $i\to\i$. Indeed, since $G(f(i),\p)$ contains no edges
joining vertices which are at a distance larger than $f(i-1)$,
the probability that a vertex $3kf(i-1)$ is a cutpoint for 
some $k=1,2,\dots,f(i)/4f(i-1)$ is larger than
$[1-\max\{a(j):j\le i-1\}] 
^{(f(i-1))^2}$ and all such events are independent. 
Hence, from (5), the number of cutpoints in $G(f(i),\p)$ is
bounded from below by the binomially distributed random variable
with  expectation ~$i$.

On the other  hand, the probability $\P(n(i),\p,\neg\psi(r))$, where
$n(i)=f(i)+i/a(i)$,  is bounded from  below by some constant
independent of $i$, which quickly tends to 1 as $r$ grows.
Indeed, call an edge
{\sl $k$-small} [{\sl $k$-large\/}]
if it is of the type $\{j,j+f(k)\}$
[$\{n(i)-j, n(i)-j-f(k)\}$]  for some
$j\le k/a(k)$.
Since $f(k)\ge (k+1)/a(k+1)$,  the existence of at least one
$k$-small and $k$-large edge for each $k=r,r+1, \dots, i$, implies
$\neg\psi(r)$. The probability  that none of the $k$-small 
[$k$-large]
edges appears in $G(n(i),\p)$ equals 
$(1-a(k))^{k/a(k)}$ and so the probability that it happens
for some $k=r,r+1, \dots,i$ is bounded from above by
$$2\sum_{k=r}^{i} (1-a(k))^{k/a(k)}
\le 2\sum_{k=r}^{\i}\exp(-k)<4e^{-r}\,.
$$
To complete the proof of Theorem 3 it is enough to observe that 
when the sequence $\a$ contains a finite number of ones we may
ignore them and repeat the above argument, whereas in the case
when in $\a$ an infinite number of ones appear one may just 
consider the property that vertex  ~1 is adjacent to $n$.\qed

Typically, when for some probabilistic model of a finite 
structure and sentence $\psi$ from language $\bar L$
convergence does not hold, it is possible to find
another sentence  $\phi$ in $\bar L$ such that the probability
of $\phi$ has  both 0 and ~1 as the limiting points.
Our next result says  that it is not the case in  $G(n,\p)$, 
provided  $\prod_{i=1}^{\infty}(1-p(i))>0$.

\proclaim{Theorem 4}
 Let $\prod_{i=1}^{\infty}(1-p(i))>0$
and let $\{\psi_{\alpha}\}_{\alpha\in A}$
be a finite set of sentences from $\LLL$.
Then there exist a subset $A'$ of $A$,  a  positive constant $\ep>0$
and a natural number $N$ such that for all $n>N$ the probability
that  $G(n,\p)$ has the property 
$\forall_{\alpha\in A'}\psi_{\alpha}
\land \forall_{\alpha\not\in A'}\neg\psi_{\alpha}$ 
is larger than $\ep$.
\endproclaim

\demo{Proof}  Let $k$ bound from above  quantifier depth of
sentences from
$A$ and ``$\pplus$'' be the operation in a family of graphs
defined in Example 1iii.  Fact 4 guarantees the existence of
an admissible graph $G$,  such that $\Th_k(G\pplus H\pplus
G)=\Th_k(G)$ for every admissible $H$
(let us recall that all graphs are treated here as $\LLL$-models).
One can easily see that  if $\prod_{i=1}^\i(1-p(i))>0$ then
 the probability that  $G(n,\p)=G\pplus H\pplus G$ for some $H$ is
bounded from below by some positive constant independent of $n$.
Thus,  the assertion follows with  $A'=A\cap\Th_k(G)$.\qed

\proclaim{Corollary} If $\prod_{i=1}^\i(1-p(i))>0$ 
then for every  $\psi$ from $\LLL$ 
$$\limsu-\limin<1.\qed$$
\endproclaim

In order to get the convergence of $\P(n,\p;\psi)$ the 
condition $\prod_{i=1}^\i(1-p(i))>0$ must be replaced by 
a significantly stronger one. 

\proclaim{Theorem 5} If $\prod_{i=1}^\i(1-p(i))^i>0$ 
then for every 
$\psi$ from $\LLL$ the probability $\P(n,p;\psi)$ converges.
\endproclaim

\demo{Proof}
Let $\psi$ be  a sentence from $\LLL$ of quantifier depth $k$
and $\p$ be a sequence for which $\prod_{i=1}^\i(1-p(i))^i>0$.
We shall  show that the sequence $\{\P(n,\p;\psi)\}_{n=1}^\i$ 
is  Cauchy.

Let  $\{\hat G(n,\p)\}_{n=1}^\i$ be a Markov process 
such that $\hat G(n,p)$ is a graph  with 
 vertex set $[n]=\{1,2,\dots,n\}$,
and for $n\ge 2$ graph $\hat G(n+1, \p)$ is such that 
\roster
\item"(i)"
if $1\le v\le w< \lfloor n/2\rfloor$ then the pair $\{v,w\}$ is an 
edge of $\hat G(n+1,p)$ if and only if  $\{v,w\}$ is an edge
of  $\hat G(n,p)$;
\item"(ii)"
 if $\lfloor n/2\rfloor< v\le w\le n+1$ then the pair $\{v,w\}$ is an 
edge of $\hat G(n+1,p)$ if and only if  $\{v-1,w-1\}$ is an edge
of  $\hat G(n,p)$;
\item"(iii)"
 if  $1\le v\le \lfloor n/2\rfloor\le w\le n+1$ and $v\neq w$ then 
$\{v,w\}$ is an edge of  $\hat G(n+1,p)$ with probability 
$p(|v-w|)$, independently for each such pair.
\endroster
Thus, roughly speaking, graph $\hat G(n+1,\p)$ is obtained from 
$\hat G(n,\p)$ by adding a new vertex in the middle of the set 
$[n]$. Clearly, we may (and will) 
identify $\hat G(n,\p)$ with $G(n,\p)$.

Now let $G$ be an admissible graph such that for every 
other admissible $H$ we have $\Th_k(G)=\Th_k(G\pplus H\pplus G)$, 
where $\pplus$ is the operation defined in the Example 1iii. 
We show first that the probability that for some $H_1$, $H_2$, $H_3$ 
$$\hat G(n,\p)=H_1\pplus G\pplus H_2\pplus G\pplus H_3\tag 6$$  
tends to 1 as $n\to\i$.

Let $k$ denote the number of vertices in $G$ and $p(G)$ be the probability 
that $G=G(k,\p)$. Moreover set $l=l(n)=\lceil \log n\rceil$ and 
for $i=0,1,\dots, l-1$, let $X_i$ be a random variable  equal to 1 when 
$$\hat G(n,\p)=H'\pplus G\pplus H''\tag 7$$
for some $H'$ with $il$ vertices and  0 otherwise. 
Then, for the expectation of $X_i$ we have 
$$\E X_i= p(G)\prod_{s=1}^{n-1}(1-p(s))^{\min\{s,k+il,n-s\}}
\prod_{s=1}^{k+il-1}(1-p(s))^{\min\{s,k,k+il-s\}}>0$$
and for $1\le i<j\le l$,
$$\align
\E X_iX_j&= \E X_i\E X_j\prod_{s=(j-i)l-k+1}^{n-1}
(1-p(s))^{-\min\{s-(j-i)l+k,k+il,n-s\}}\\
&\le\E X_i\E X_j\prod_{s=l}^{\i}(1-p(s))^{-s},
\endalign$$
where here and below we assume that $n$ is large enough to have  $l> k$.
Thus,  the expectation of the random variable 
  $X=\sum_{i=0}^{l-1} X_i$ is of the order $\log n$ and
$$\Var X\le (\E X)^2 \(\prod_{s=l}^\i(1-p(s))\super{- s}-1\)
+O(\log n)\;, $$
so, due to Chebyshev's inequality, 
(7) holds for some $H'$ with at most $l(l-1)$ vertices 
with probability  at least 
$1-O(1/\log n)\allowbreak-
O\(\prod_{s=l}^\i(1-p(s))^{-s}-1\)$.
Clearly, an analogous argument shows that (7) remains
valid for  $H''$  of size not larger than $l(l-1)$ so 
with probability at least 
$1-O(1/\log n)-O\(\prod_{s=l}^\i(1-p(s))\super{- s}-1\)$
(6) holds for some  $H_1$ and $H_3$, both of them 
with not more  than $l(l-1)$ vertices.
   
Now assume that (6) is valid and let $m>n$.  Then, 
the probability that for some  $H'_2$ we have 
$$\hat G(m,\p)=H_1\pplus G\pplus H'_2\pplus G\pplus H_3\tag{$6'$}$$  
with the same $H_1$ and $H_3$ as in (6) is at least
$$1-\(\prod_{s=n/2-l^2}^{m-1}
(1-p(s))^{\min\{s-n/2+l^2,l^2,m-s\}}\)^2
\ge 1-\prod_{s=n/3}^\i(1-p(s))^{2s}\;.$$
But, provided (6) and ($6'$) holds, 
$\Th_k(\hat G(n,\p))=\Th_k(\hat G(m,\p))$. 
Hence, for every $n$ and $m$ such that $m>n$ we have 
$$|\P(n,\p,\psi)-\P(m,\p,\psi)|  \le \ep(n)\;,$$
where 
$$\ep(n)=O(1/\log n)+O\(\prod_{s=(\lceil\log n\rceil)^2}
^\i(1-p(s))\super{- s}-1\)+ 
       1-\prod_{s=n/3}^\i(1-p(s))^{2s}\to 0\;.$$
Thus, sequence $\{\P(n,\p,\psi)\}_{n=1}^\i$, being Cauchy,
must converge.\qed

\subheading{6. First order properties of $C(n,\p)$}

Let $C(n,\p)$  denote a graph with  vertex set
$[n] $ in which a pair of vertices $v$, $w$, are joined
 by an edge with probability
$p(\min\{|v-w|, n-|v-w|\})$. 
Let $L^c$ be the first order logic which uses
only the adjacency predicate, in $L^c_+$ one may  say  also
 $v=w+1$\,(mod$(n)$), whereas the vocabulary of $L^c_\le$
contains the predicate
$C(v_1,v_2,v_3)$ which means that starting from ~$v_1$ and
moving clockwise $v_2$ is met before $v_3$, 
i.e. for some cyclic permutation $\sigma$
of indices $v_{\sigma(1)}\le v_{\sigma(2)}\le v_{\sigma(3)}$.
 It turns out that differences between $L^c$, 
 $\LL^c$ and $\LLL^c$
are not so substantial as those between $L$, $\LL$ and $\LLL$.

\proclaim{Theorem 6}
\roster
\item"(i)"
If a sequence $\p$ contains only finitely many non-zero terms
then for every sentence $\psi$ from $\LLL^c$ a zero-one  law holds.
\item"(ii)"
For every infinite sequence $\a$ such that $0<a(i)<1$, 
$i=1,2,\dots$,  there exist a sequence $\p$ obtained from $\a$ 
by the addition of some
number of zero terms and a  first order sentence $\psi$
from $L^c$ such that $\limin=0$ but $\limsu=1$.
\endroster
\endproclaim

\demo{Proof}
The first part of Theorem 6 follows  from the Lemma in a similar way
as in the proof of Theorem~3i.
 To show (ii) assume, for simplicity, that the sequence $a(i)$
decreases,   define  $\p$ setting
$$p(j)=\cases a(i)&\quad\text{ if } j=\lfloor 3\super
i/a(i)\rfloor\\   0& \quad \text { otherwise}\,,
\endcases$$
and let $\psi$ be a sentence that $C(n,\p)$ contains a cycle of
length 3.  
It is not hard to see that $C(3\lfloor 3^{i}/a(i)\rfloor-1,\p)$
contains no cycles of length 3 whereas the number of such  
cycles in  $C(3\lfloor 3\super{i}/a(i)\rfloor,\p)$
is binomially distributed with parameters 
$\lfloor 3^{i}/a(i)\rfloor $ and $a(i)$.\qed

Clearly, the proof of Theorem 6ii is based on the fact that, unlike in
the case of $G(n,\p)$,   for
some subgraphs $H$ the probability that $H$ is contained in $C(n,\p)$
might be smaller than the probability that $H$ is contained in
$C(n+1,\p)$.
To eliminate   such a pathological situation let us
call a subgraph  $H$  of $C(n,\p)$ with vertices
$v_1,v_2,\dots,v_k$ $L^c$-{\sl flat} [$\LL^c$-{\sl flat\/}]
 if there is a sequence $w_1,w_2,\dots,w_k$ of vertices of
$G(n,\p)$ such that for every $i,j\in [k] $
the probability $p(|w_j-w_i|)$ is positive
whenever $\{v_i,v_j\}$ is an edge of $H$
[and, moreover, $w_j=w_i+1$  if and only if
$v_j=v_i+1$\,(mod$(n)$)]. Furthermore,  such a
subgraph $H$
is $\LLL^c$-{\sl flat} if there exists $i$, $1\le  i \le k$,
such that for  all  $i'$, $i''$, $i'\neq i''$, 
$1\le i',i''\le k$, such that   $C(v_i,v_{i'},v_{i''})$ and 
$ |v_{i'}-v_{i}|+|v_{i}-v_{i''}|\le n/2$, the  pair
$\{v_{i'},v_{i''}\}$ is not an edge of $H$.
Finally call a sequence $\p$ {\sl flat} 
[{\sl asymptotically flat}] 
with respect to language $L_\bullet$,
$L_\bullet=L^c, \LL^c,\LLL^c$, if
the probability that  $C(n,\p)$ contains a subgraph which is not
$L_{\bullet}$-flat is 0 [tends to 0].

\proclaim{Theorem 7}
If a sequence $\p$ fulfilling the assertion of Theorem 1i is
asymptotically flat with respect to $L_\bullet$, 
where  $L_\bullet=L^c,
\LL^c$ or $\LLL^c$,
then for every $\psi$ from  $L_\bullet$ a zero-one law holds.
\endproclaim

\demo{Proof}
 Let $\psi$  be a sentence of quantifier depth $k$.
From Fact 4, there exists a  $L_\bullet$-flat graph $G$ such that
for every  $L_\bullet$-flat $H$ we have  
$\Th_k(G)=\Th_k(G\pplus H)$,where both graphs $G$ and $G\pplus H$
are treated  as  $L_\bullet$-models.
Now it is enough to observe that from the Lemma the probability that
$C(n,\p)$ contains an exact copy of $G$ tends to 1 as $n\to\i$,
provided  the assertion of Theorem 1 holds.\qed

\proclaim{Corollary} If for each $k$ there exists $m$ such that
$p(im)>0$ for $i=1,2,\dots,k$, and if the sequence $\p$ fulfills
 condition (2) then for every $\psi$ from $\LL^c$ a zero-one
law holds.
\endproclaim

\demo{Proof} It is enough to note that each sequence $\p$ for
which the  assumption of the Corollary remains valid is 
$\LL^c$-flat.\qed

\subheading{7. Final remarks and comments}

 One may ask whether additional restrictions imposed on the
sequence $\p$ like non-negativity of all terms or monotonicity
could lead to other interesting
results concerning convergence properties of $G(n,\p)$.
However, all sequences $\p$ which appeared in all our
 counterexamples (with the single exception  of  Theorem 6ii,
where non-negativity plays an important role)
could be modified in such a way that they become both 
non-negative and monotonically decreasing, so no new 
sufficient conditions for convergence
can be shown under these new assumptions.

In the paper we have studied  properties of a random graph  $G(n,\p)$
which is a generalization  of $G(n,p)$ when the probability $p$ does not
depend on $n$.
The  problem of characterizing convergence properties in the case
when $\p$ varies with $n$ seems to be a  much more challenging
problem -- we recall only that if  $p(n)\to 0$ 
then convergence properties of $G(n,p)$ become quite involved 
and strongly  depend on the limit behaviour
of function $p=p(n)$ (for details see papers of Shelah and Spencer
[SS 88] and \L uczak and Spencer [\L S 91]).

It is not hard to observe that  if we are interested only 
in properties described by sentences of quantifier depth
bounded by $k$ then Theorems 1 and 2 remain valid also when 
zeros are replaced by very small constants $\ep(k)$
i.e. (2) could be replaced by
$$-\log\Big(\prod_{i=1}^n(1-p(i)\Big)\big/(\log
n)<\ep\,.\tag{$2'$}$$
On the other hand, even if $\p$ is such that the probability
$\P(n,\p;\psi)$ converges for every sentence $\psi$ from $L$ the
rate of this convergence may be very slow for sentences of
large quantifier depth.

\proclaim{Theorem 8}
There exists a sequence $\{\psi_k\}_{k=1}^\i$ of
first order sentences from $L$,  $\psi_k$ of depth $k$
for every $k=1,2,\dots$, such that
for every sequence $\p$, $0<p(i)<1$, for which  (2) holds we have
$$\lim_{n\to\i}\P(n,\p;\psi_k)=1\quad\quad
\text{for every }\quad k=1,2,\dots$$  
but the function $m(k)=\min\{i:\P(i,\p;\psi_k)>0\}$
grows faster than any recursive function of $k$.
\endproclaim

\demo{Proof}
It is well known (see, for instance, [Tr 50]) 
that there exists a sequence $\{\phi_k\}_{k=2}^\i$
such that for $k=2,3,\dots$,  $\phi_k$ is a first order sentence
of  depth $k$ from $L$ and a function
$m'(k)=\min\{\operatorname{card}(M)\,:\, 
M\text{ is  a model of $\phi_{k-1}$}\}$
grows faster than any recursive function of $k$. Let $\psi_k$
be the sentence which states that a graph contains a vertex $v$
such that for a subgraph induced by all neighbours of $v$
$\phi_{k-1}$ holds. Clearly, $\psi_k$ has depth $k$. Now let $\p$
be a sequence such that $0<p(i)<1$ for all $i=1,2,\dots$ for which
 (2) holds. For such a sequence all graphs are admissible, so
 the Lemma implies that  $\P(n,\p;\psi_k)$ tends to 1 as $n\to\i$.
On the other hand,
$$m(k)=\min\{i:\P(i,\p;\psi_k)>0\}\ge m'(k)\,.\qed$$

Thus, the behaviour of $ \P(n,\p;\psi)$ for small $n$ does not tell us 
very much about the asymptotic behaviour of $ \P(n,\p;\psi)$. 
It is not hard to show that even if $\p$ is such that 
for every first order sentence a zero-one law holds, 
there is no procedure  which,  applied to $\psi$, 
could decide whether $\lim_{n\to\i}\P(n,\p;\psi)=0$ or  
$\lim_{n\to\i}\P(n,\p;\psi)=1$.    

\proclaim{Theorem 9}
Let $\p$ be a sequence such that  $0<p(i)<1$ for all $i$, 
for which (2) holds. Then there exists no procedure which 
can decide for each  first order sentence $\psi $ from $L$, 
whether  $\lim_{n\to\i}\P(n,\p;\psi)=0$ or  
$\lim_{n\to\i}\P(n,\p;\psi)=1$.    
\endproclaim

\demo{Proof} 
Let $\phi$ be a first order sentence from $L$ and $\psi_\phi$ 
denote the sentence
that for some vertex $v$ in a graph the subgraph which is induced 
in a graph by all neighbours of $v$ has property $\phi$.
Since $p(i)>0$ for all $i$, every graph 
is admissible for $G(n,\p)$  and the Lemma implies that,
with probability tending to 1 as $n\to\i$, 
every finite graph appears in $G(n,\p)$ as a component.
Thus,  $\lim _{n\to\i} \P(n,\p;\psi_\phi)=1$  if and only if 
$\phi$ is satisfied for some finite graph. 
Now the assertion follows from the fact that, due to 
the Traktenbrot-Vought Theorem [Tr 50], there is  no decision procedure
to determine whether a first order sentence 
$\phi$ from $L$ has a finite model.\qed

\medskip
{\bf Acknowledgement.}\, We wish to thank anonymous referees 
for their numerous insightful remarks 
on the first version of the paper.


\enddocument